\documentclass[12pt,a4paper]{amsart}
\usepackage[utf8]{inputenc}
\usepackage{amsfonts}
\usepackage{amssymb}
\usepackage{bbm}
\usepackage{comment}
\usepackage{enumitem}
\usepackage{dsfont}
\usepackage{amsthm}
\usepackage{amsmath}
\usepackage{commath}
\usepackage{url}
\usepackage{epsfig}
\usepackage{adjustbox}
\usepackage[ruled,vlined,linesnumbered]{algorithm2e}
\usepackage{cool}
\usepackage[table]{xcolor}
\usepackage{tabularx}
\usepackage[a4paper]{geometry}
\usepackage{caption}
\usepackage{booktabs}
\usepackage{subcaption}
\usepackage{graphicx}
\usepackage{a4wide}
\usepackage{todonotes}
\usepackage{scrextend}
\usepackage{mathtools}
\mathtoolsset{showonlyrefs,showmanualtags}

\newcommand{\ls}{\leqslant}
\newcommand{\gs}{\geqslant}
\newtheorem{thm}{Theorem}
\newtheorem{corollary}[thm]{Corollary}
 \newtheorem{proposition}[thm]{Proposition}
\newtheorem{lem}[thm]{Lemma}

\title{On extreme values for the Sudler product of quadratic irrationals}
\author{Manuel Hauke}
\address{TU Graz, Austria}
\email{hauke@math.tugraz.at}
\date{\today}

\begin{document}

\maketitle

\begin{abstract}
     Given a real number $\alpha$ and a natural number $N$, the Sudler product is defined by
     \mbox{$P_N(\alpha) = \prod_{r=1}^{N} 2 \left\lvert \sin(\pi\left(r\alpha \right))\right\rvert.$}
     Denoting by $F_n$ the $n$--th Fibonacci number and by $\phi$ the Golden Ratio, we show that for $F_{n-1} \ls N < F_n$, we have
     $P_{F_{n-1}}(\phi)\ls P_N(\phi) \ls P_{F_{n}-1}(\phi)$ and $\min_{N \gs 1} P_N(\phi) = P_1(\phi)$, thereby proving a conjecture of Grepstad, Kaltenböck and Neumüller. Furthermore, we find closed expressions
     for $\liminf_{N \to \infty} P_N(\phi)$ and $\limsup_{N \to \infty} \frac{P_N(\phi)}{N}$ whose numerical values can be approximated arbitrarily well. We generalize these results to the case of quadratic irrationals $\beta$ with continued fraction expansion $\beta = [0;b,b,b\ldots]$ where $1 \ls b \ls 5$,
     completing the calculation of $\liminf_{N \to \infty} P_N(\beta)$, $\limsup_{N \to \infty} \frac{P_N(\beta)}{N}$ for $\beta$ being an arbitrary quadratic irrational with continued fraction expansion of period length 1.
\end{abstract}

\section{Introduction and statement of results}
For $\alpha \in \mathbb{R}$ and $N$ a natural number, the Sudler product is defined as
\begin{equation}\label{sudler}P_N(\alpha) := \prod_{r=1}^{N} 2 \left\lvert \sin \pi r\alpha \right\rvert.
\end{equation}
This product first appeared in a paper of Erd\"os and Szekeres \cite{erdos_szekeres}, where 
it was proven that

\[\liminf_{N \to \infty} P_N(\alpha) = 0, \quad \limsup_{N \to \infty} P_N(\alpha) = \infty\]
holds for almost every $\alpha$.
In the same work, they raised the question whether
\begin{equation}\label{erdos_conjecture}\liminf_{N \to \infty} P_N(\alpha) = 0\end{equation}
holds for all $\alpha \in \mathbb{R}$. Note that $\limsup\limits_{N \to \infty} P_N(\alpha) = \infty$ cannot hold for all $\alpha$ since for rational $\alpha = \frac{n}{m}$,
$P_N(\alpha) = 0$ for $N \gs m$. Thus, by periodicity of the sine function, we can restrict the asymptotic analysis of $P_N(\alpha)$ to the case of 
irrational numbers $\alpha$ in the unit interval.\vspace{5mm}

Denoting by $\| P_N \| = \max_{0 < \alpha < 1} P_N(\alpha)$, Erdös and Szekeres \cite{erdos_szekeres} claimed that the limit
$\lim_{N \to \infty} \| P_N \|^{1/N}$ exists and equals a value between $1$ and $2$, without formally proving it. This was done by Sudler \cite{sudler} and Wright \cite{wright} showing that 
$\lim_{N \to \infty} \| P_N \|^{1/N} = C \approx 1.22$.
Inspired by this, the order of growth of Sudler products was extensively examined from a metric point of view. For more results in this area, we refer
the reader to \cite{atkinson,bbr,bell,bc,fh,kol,kol2}.\vspace{5mm}

In this paper, we consider the pointwise behaviour of Sudler products.
In contrast to the result on $\lim_{N \to \infty} \| P_N \|^{1/N}$, Lubinsky and Saff \cite{ls} showed that for almost every $\alpha$, \mbox{$\lim_{N \to \infty} P_N(\alpha)^{1/N} = 1$}.
Estimates on $P_N(\alpha)$ for fixed $\alpha$ were used by Avila and Jitomirskaya \cite{aj}
to solve the Ten Martini Problem, and also play a role in a proof of Avila, Jitomirskaya and Marx in \cite{ajm}.
Recently, Aistleitner and Borda \cite{quantum_invariants} established a connection between Sudler products and the work of Bettin and Drappeau \cite{bd1,bd2,bd3} on the order of magnitude of the Kashaev invariant of certain hyperbolic knots following the work of Zagier \cite{zag}.
\vspace{5mm}

 Mestel and Verschueren \cite{mestel} examined the behaviour of $P_N(\phi)$ where $\phi = \frac{\sqrt{5}-1}{2} = [0;1,1,1,\ldots]$ is the fractional part of the Golden Ratio.
Throughout this paper, $\phi$ will always denote this value and $P_N(\phi)$ will be called the Sudler product of the Golden Ratio. %
 The authors of \cite{mestel} showed that the limit along the Fibonacci sequence $\lim\limits_{n \to \infty} P_{F_n}(\phi)$ exists, without giving a closed expression of its value. %
 Additionally, it was proven in \cite{mestel} that \begin{equation}\label{phi_n}\lim\limits_{n \to \infty} \frac{P_{F_n-1}(\phi)}{F_n-1} = \frac{\sqrt{5}}{2\pi}\lim\limits_{n \to \infty} P_{F_n}(\phi).\end{equation}
 Note that the Fibonacci numbers are the denominators of the continued fraction convergents of $\phi$, hinting at a connection between the Sudler product of $\alpha$ and its Diophantine approximation properties. This was established more broadly by Aistleitner, Technau and Zafeiropoulos \cite{tech_zaf} who generalized Mestel and Verschueren's work to quadratic irrationals of the form $\beta(b) = [0;b,b,\ldots]$ for arbitrary natural numbers $b$. They showed that 
 for $q_n$ being the denominator of the $n$--th continued fraction convergent,
 the limit
 $ \lim_{n \to \infty} P_{q_n}(\beta) = C_b > 0$ exists, also giving a closed expression for $C_b$. Further, they generalized \eqref{phi_n} by showing that
 \begin{equation}\label{limsup_mirror}\lim\limits_{n \to \infty} \frac{P_{q_n-1}(\phi)}{q_n-1} = \frac{\sqrt{b^2+4}}{2\pi}C_b.\end{equation}
In fact $q_n = q_n(\beta)$ depends on $\beta$, but we will suppress this dependence during the rest of this paper for the sake of readability.\vspace{5mm}\par

Returning to \eqref{erdos_conjecture}, Lubisky \cite{lubinsky} showed that 
$\liminf_{N \to \infty} P_N(\alpha) = 0$ holds for any $\alpha \in \mathbb{R}$ that has unbounded 
continued fraction coefficients, conjecturing it to hold for any $\alpha$.
However, this conjecture was disproven by Grepstad, Kaltenböck and Neumüller \cite{grepstad_neum} as they 
proved that
\begin{equation}\label{positive_liminf}\liminf_{N \to \infty} P_N(\phi) > 0.\end{equation}
Later Aistleitner, Technau and Zafeiropoulos \cite{tech_zaf} found a very close connection between the behaviour of $\liminf\limits_{N \to \infty} P_N(\alpha)$ and $\limsup\limits_{N \to \infty} \frac{P_N(\alpha)}{N}$. With this connection, they extended the counterexample from \cite{grepstad_neum}
 to all irrationals of the form $\beta(b) = [0;b,b,\ldots]$ with $1 \ls b \ls 5$, showing that there is a sharp threshold for the value $b$ where the behaviour of $\liminf\limits_{N \to \infty} P_N(\beta)$ and $\limsup\limits_{N \to \infty} \frac{P_N(\beta)}{N}$ changes.\vspace{5mm}
 
\noindent
{\bf Theorem A }(Aistleitner, Technau and Zafeiropoulos \cite[Theorem 6]{tech_zaf}).
 Let $b$ be a positive integer and let $\beta = \beta(b) = [0;b,b,\ldots]$. Then the following holds.

 \begin{enumerate}
     \item[(i)] If $b \ls 5$, then 
     $\liminf\limits_{N \to \infty} P_N(\beta) > 0$ and $\limsup\limits_{N \to \infty} \frac{P_N(\beta)}{N} < \infty$.
     \item[(ii)] If $b \gs 6$, then 
     $\liminf\limits_{N \to \infty} P_N(\beta) = 0$ and $\limsup\limits_{N \to \infty} \frac{P_N(\beta)}{N} = \infty$.
 \end{enumerate}
  However, the precise values of 
 $\liminf\limits_{N \to \infty} P_N(\beta),\; \limsup_{N \to \infty} \frac{P_N(\beta)}{N}$
 for $1 \ls b \ls 5$
 have not been determined so far. They will be established by the following theorem.
 
 \begin{thm}\label{main_3}
Let $b \ls 5$ be a positive integer and let $\beta = \beta(b) = [0;b,b,\ldots]$. Then we have
     \begin{align}\label{thm_3_eq}
     \liminf_{N \to \infty} P_N(\beta) &= \lim_{n \to \infty} P_{q_n}(\beta) = C_b, \quad
     \\\limsup_{N \to \infty} \frac{P_{N}(\beta)}{N} &= \lim_{n \to \infty} \frac{P_{q_n-1}(\beta)}{q_n-1} = \frac{\sqrt{b^2 + 4}}{2\pi}C_b,\end{align}
     where 
    \begin{equation}\label{C_b}
    C_b  = \frac{2\pi}{\sqrt{b^2 + 4}}\prod_{n=1}^{\infty}
    \left(1 - \frac{1}{\sqrt{b^2 + 4}}\frac{\left\{n\beta\right\} - \frac{1}{2}}{n}\right)^2 - \frac{\left(\frac{1}{2\sqrt{b^2 +4}}\right)^2}{n^2}.
    \end{equation}
\end{thm}
It was already shown in \cite{tech_zaf} that $\lim_{n \to \infty} P_{q_n}(\beta) = C_b$, although 
the closed expression differs from \eqref{C_b}. The expression used in Theorem \ref{main_3} first appeared in a more general setting in \cite{quantum_invariants} where arbitrary quadratic irrationals were considered.

\subsection{Extreme values in a finite range}

We proceed to discuss the behaviour of $P_N(\beta)$ not only asymptotically, but also for finite intervals for $N$.
  In their survey paper on Sudler products \cite{grepstad_survey}, Grepstad and Neumüller 
  looked at $P_N(\phi)$ for $N$ between $1$ and a reasonable large $N_0$.
  \begin{figure}[h]
  \centering
  \includegraphics[width=.7\linewidth]{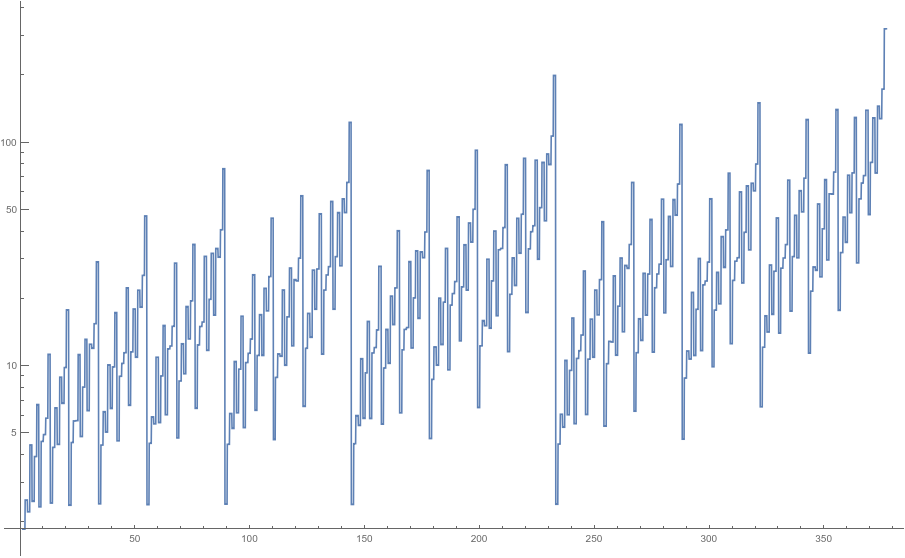}
  \caption{Values of $P_N(\phi)$ for $1 \ls N \ls F_{14}-1 = 376$ on a logarithmic scale.}
\end{figure}
They observed that the peaks of $P_N(\phi)$ seem to be at $N = F_{n}-1$ with a local minimum
at $N = F_n$ for every $n \in \mathbb{N}$. Additionally, they observed that $P_N(\phi)$ appears to be minimal when $N =1$,
an observation that extends to $\beta(b)$ with $1 \ls b \ls 5$.
This led them to the following two conjectures.
    \begin{enumerate}
        \item[(i)]
  Let $n \gs 3, F_{n-1} \ls N < F_n$ where $F_n$ denotes the $n$--th Fibonacci number. Then we have
  \begin{equation}
     \label{conj_grep_1}
 P_{F_{n-1}}(\phi)\ls P_N(\phi) \ls P_{F_{n}-1}(\phi).
 \end{equation}
    \item[(ii)] %
  Let $1 \ls b \ls 5$, $\beta = \beta(b) = [0;b,b,\ldots]$.
  Then we have

 \begin{equation}\label{thm_2_eq}\min_{N \gs 1} P_N(\beta) = P_1(\beta).\end{equation}
 \end{enumerate}
We will prove these conjectures here, also extending \eqref{conj_grep_1} to $\beta(b)$ with $1 \ls b \ls 5$.\vspace{5mm}

\begin{thm}\label{main_1}
 
 Let $1 \ls b \ls 5$, $\beta = \beta(b) = [0;b,b,\ldots]$.
 The following statements hold.
 
 \begin{enumerate}
     \item[(i)]
     Let $N$ be a natural number such that $q_n \ls N < q_{n+1}$, $n \gs 1$. Then we have
       \begin{equation}\label{thm_1_2}P_{q_n}(\beta) \ls P_N(\beta)\ls P_{q_{n+1}-1}(\beta).
  \end{equation}
  \item[(ii)] We have
  \begin{equation}\label{min_at1}\min_{N \gs 1} P_N(\beta) = P_1(\beta).\end{equation}
 \end{enumerate}
 \end{thm}
 \captionsetup[subfigure]{labelformat=empty}
\begin{figure}[h]
\begin{subfigure}{.5\textwidth}
  \centering
  \includegraphics[width=.9\linewidth]{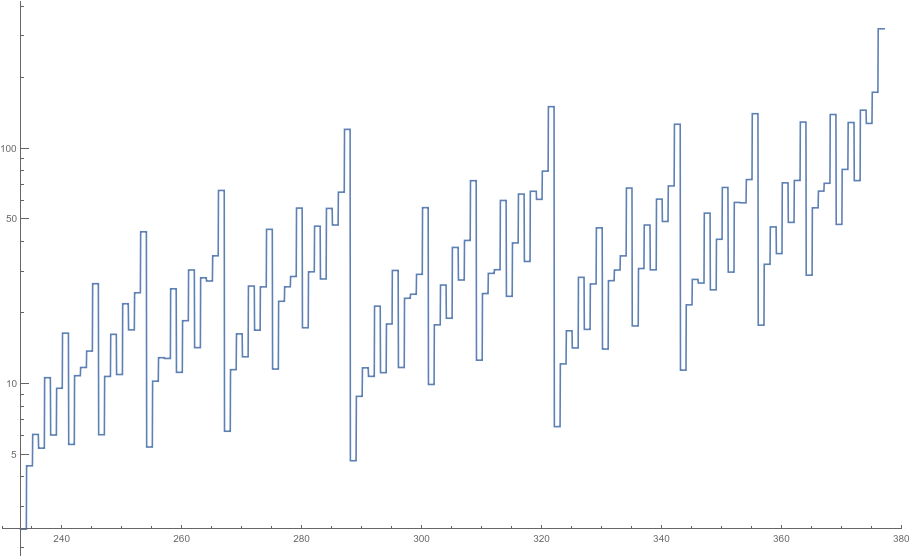}
  \caption{\phantom{ }\hspace{.7cm} Values of $P_N(\phi)$ for $F_{13} \ls N \ls F_{14}$\newline
  \phantom{ }\hspace{.7cm} on a logarithmic scale.
  }
\end{subfigure}%
\begin{subfigure}{.5\textwidth}
  \centering
  \includegraphics[width=.9\linewidth]{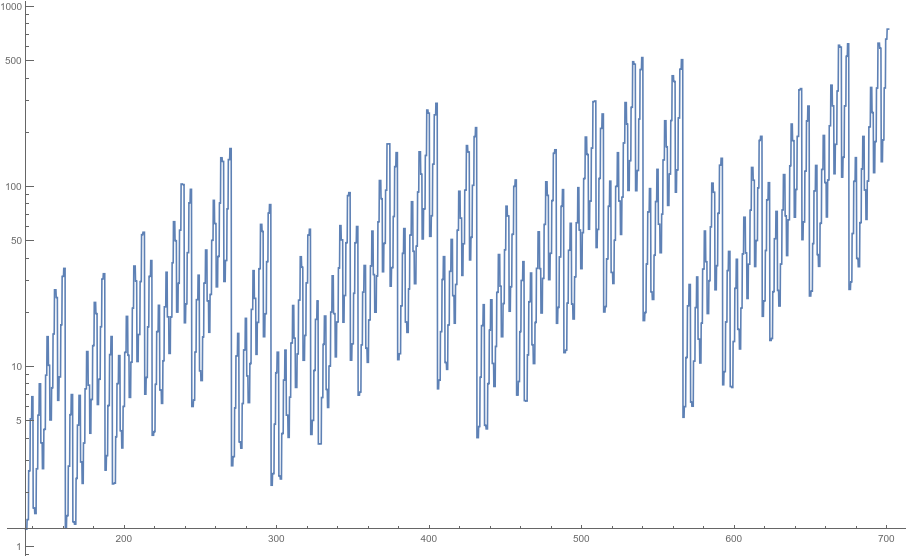}
  \caption{\phantom{ }\hspace{.7cm} Values of $P_N(\beta(5))$ for $q_3 \ls N \ls q_4$\newline
  \phantom{ }\hspace{.7cm} on a logarithmic scale.}
\end{subfigure}
\caption{Visualization of \eqref{thm_1_2} for $\phi$ and $\beta(5)$.
}
\end{figure}
 Both statements (i) and (ii) of Theorem \ref{main_1} are optimal for quadratic irrationals with period length $1$ in the sense that they cannot be extended to $\beta(b)$ with $b \gs 6$:
 By Theorem A, we know that in that case $\liminf_{N \to \infty} P_N(\beta) = 0$, already making \eqref{min_at1} impossible to hold. Since Aistleitner, Technau and Zafeiropoulus showed in \cite{tech_zaf} that for any $b \in \mathbb{N}$ we have $\lim_{n \to \infty} P_{q_n}(\beta) = C_b > 0$, also $P_{q_n}(\beta) \ls P_N(\beta)$
 must be violated for any $b \gs 6$ and some $n,N$ sufficiently large. By a similar argument, we see that
 $P_N(\phi) \ls P_{q_n-1}(\beta)$ eventually fails as well.
 
 \begin{figure}[h]
\begin{subfigure}{.5\textwidth}
  \centering
  \includegraphics[width=.9\linewidth]{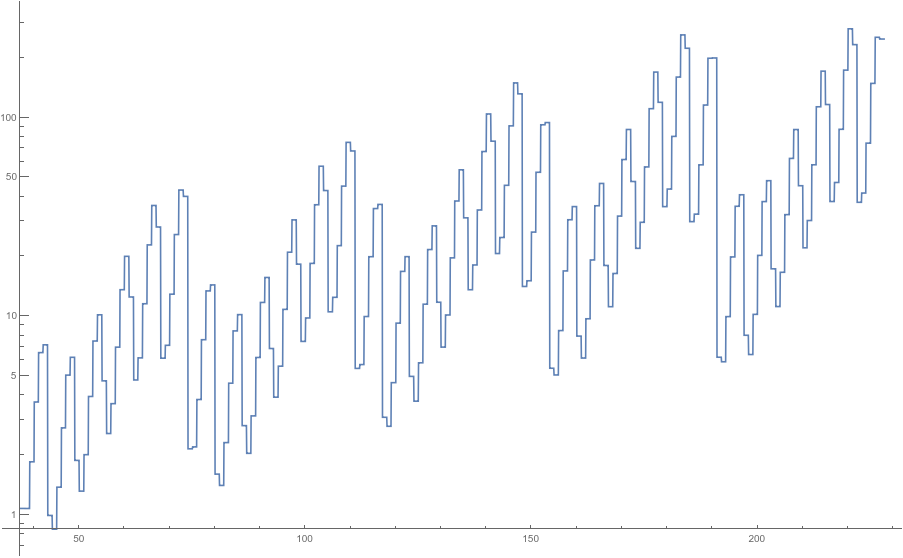}
  \caption{\phantom{ }\hspace{.7cm} Values of $P_N(\beta(6))$ for $q_2 \ls N \ls q_3-1$\newline
  \phantom{ }\hspace{.7cm} on a logarithmic scale.}
\end{subfigure}%
\begin{subfigure}{.5\textwidth}
  \centering
  \includegraphics[width=.9\linewidth]{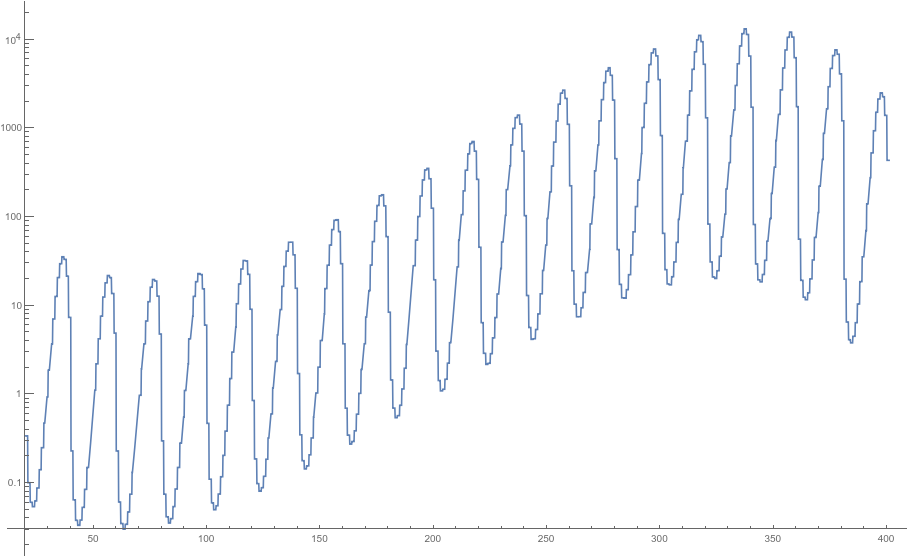}
  \caption{\phantom{ }\hspace{.7cm} Values of $P_N(\beta(20))$ for $q_1 \ls N \ls q_2-1$ \newline
  \phantom{ }\hspace{.7cm} on a logarithmic scale.}
\end{subfigure}
 \caption{We see that the peak of $P_N(\beta)$ does not appear at $q_n -1$ when $b \gs 6$.
 This is in accordance with the results of \cite{quantum_invariants} which describe the behaviour of the $P_N(\alpha)$ for irrational numbers $\alpha$ which have ``large'' partial quotients in their continued fraction expansion.}
 \end{figure}
  
  \subsection{Directions for further research}
  Recently, Grepstad, Neumüller and Zafeiropoulos \cite{grepstad_new} studied the behaviour of $\liminf_{N \to \infty} P_N(\alpha),\;  \limsup_{N \to \infty} \frac{P_N(\alpha)}{N}$ for quadratic irrationals $\alpha$ of the form $\alpha = [0;\overline{a_1,\ldots, a_{\ell}}] = 
  [0;a_1,\ldots, a_{\ell},a_1,\ldots,a_{\ell},a_1,\ldots]$. They proved that if \mbox{$a_k := \max\limits_{1 \ls i \ls \ell} a_i\gs 23$}, then 
  \begin{equation}\label{not_int_case}\liminf_{N \to \infty} P_N(\alpha) = 0, \quad \limsup_{N \to \infty} \frac{P_N(\alpha)}{N} = \infty.\end{equation}
  Actually, they show that the properties in \eqref{not_int_case} are invariant under the Gauss map\\
  \mbox{$[0;a_1,a_2,\ldots] \mapsto [0;a_2,a_3,\ldots]$}, so the statement can be extended to 
  arbitrary quadratic irrationals $\alpha = [0;a_1,\ldots, a_p,\overline{a_{p+1},\ldots, a_{p+\ell}}]$
  with $a_k := \max\limits_{p+1 \ls i \ls p+\ell} a_i$.
      It remains an open question how far the condition $a_k \ls 23$ can be relaxed. Numerical experiments suggest that $a_k \gs 6$ is sufficient to have \eqref{not_int_case}, which would be in accordance with Theorem A (ii) where the period length is $1$.
        Grepstad and Neumüller \cite{gn_quad} showed that for purely periodic quadratic irrationals with period length $\ell$, and $1 \ls r \ls \ell$, $\lim_{N \to \infty} P_{q_{m\ell + r}} = C_r > 0$ exists,
    hence we can argue similarly to before that whenever \eqref{not_int_case} is fulfilled, the statements
    from Theorem \ref{main_1} cannot hold. However, the converse does not need to be true since Theorem \ref{main_1} considers the behaviour of $P_N(\beta)$ in finite ranges. So one can ask which quadratic irrationals fulfill $\liminf_{N \to \infty} P_N(\beta) > 0$ and $\limsup_{N \to \infty} \frac{P_N(\beta)}{N} < \infty$, but do not fulfill \eqref{thm_1_2} and \eqref{min_at1}.
    It is clear that this holds for irrationals with large partial quotients in the pre-period, so a characterization of quadratic irrationals fulfilling \eqref{thm_1_2} and \eqref{min_at1} will probably
    only consider irrationals with purely periodic continued fraction or will need to consider the value of
    $a_k := \max\limits_{i \gs 1} a_i$.
    
    \subsection{Notation}\label{notation} Given two functions $f,g:(0,\infty)\to \mathbb{R},$ we write $f(x) = \mathcal{O}(g(x))$, \mbox{$f(x)= o(g(x))$} and $f(x)=\Omega(g(x))$ when  \vspace{-2mm}
\[\limsup_{x\to\infty} \frac{|f(x)|}{|g(x)|} < \infty, \quad   \quad \lim_{x\to\infty} \frac{f(x)}{g(x)} =0 \quad \text{ and } \quad \limsup_{x\to\infty}\frac{f(x)}{g(x)}>0 \] respectively. Given a real number $x\in \mathbb{R},$ we write $\{x\}$ for the fractional part of $x$ and \mbox{$\|x\|=\min\{|x-k|: k\in\mathbb{Z}\}$} for the distance of $x$ from its nearest integer.
For $A$ being a logical expression, we write
\[\mathds{1}_A = \begin{cases} 1 \text{ if } A \text{ is true} \\
0 \text{ otherwise}.\end{cases}\]

\section{Preliminaries}
\subsection{Continued fractions}
 
We will quickly recall all facts needed in this paper about continued fractions, stating several identities used in the subsequent proofs. For a more detailed background, see e.g. \cite{all_shall, rock_sz, schmidt}.
Every irrational $\alpha$ has a unique infinite continued fraction expansion $[a_0;a_1,...]$ with  convergents $p_k/q_k = [a_0;a_1,...,a_k]$ fulfilling the recursions 
\begin{equation}\label{recursions}p_{k+1} = p_{k+1}(\alpha) = a_{k+1}p_k + p_{k-1}, \quad q_{k+1} = q_{k+1}(\alpha) = a_{k+1}q_k + q_{k-1}\end{equation}
with initial values $p_0 = a_0,\; p_1 = a_1a_0 +1,\; q_0 = 1,\; q_1 = a_1$.
From \eqref{recursions}, one can deduce that
$q_kp_{k-1} - p_kq_{k-1} = (-1)^k$ holds for any $k \gs 1$, so in particular,
\begin{equation}\label{conv_ident}p_kq_{k-1} \equiv (-1)^{k+1} \mod{q_k}.\end{equation}
Furthermore, $q_k,q_{k-1}$ are coprime and thus,
\begin{align}\label{bijection}
    g: \{1,\ldots, q_{k}-1\} &\to \{1,\ldots, q_{k}-1\}
    \\n &\mapsto q_{k-1}n \mod {q_k}
\end{align}
is bijective. Additionally, we have
\begin{equation}\label{alpha_deltak} \alpha = \frac{p_k}{q_k} + (-1)^k\frac{\|q_k\alpha\|}{q_k}.\end{equation}
We know that $p_k/q_k$ approximates $\alpha$ very well. Indeed, we have the following well-known inequalities for $k \gs 1$:

\begin{equation}\label{qkqkalpha}
    \frac{1}{a_{k+1}+2}\ls q_k\|q_k\alpha\| \ls \frac{q_k}{q_{k+1}} \ls \frac{1}{a_{k+1}},
\end{equation}
\begin{equation}\label{best_approx}
    \|q\alpha\| > \|q_k\alpha\|, \quad 1 < q < q_{k+1}, q \neq q_k,
\end{equation}

\begin{equation}\label{approx_quality}
\Big\lvert \frac{p_k}{q_{k}} - \beta\Big\rvert \ls\frac{1}{q_k^2}.
\end{equation}
\vspace{5mm}

Fixing an irrational $\alpha = [a_0;a_1,...]$, the Ostrowski expansion of a non-negative integer $N$ is the unique representation

\begin{equation}\label{ostrowski}N = \sum_{\ell = 0}^k b_{\ell+1}q_{\ell} \quad \text{ where }
b_{k+1} \neq 0, 0 \ls b_1 < a_1, \quad 0 \ls b_{\ell} \ls a_{l} \text { for } \ell \gs 2,
\end{equation}
with the additional rule that
$b_{\ell-1} = 0$ whenever $b_{\ell} = a_{\ell+1}$.\vspace{5mm}\par

We proceed to the special case of irrationals of the form $\beta = \beta(b) = [0;b,b,\ldots]$. Here we have the following well-known identities:

\begin{align}
    \label{beta_identity_0}
    &\beta = \frac{1}{2}(-b + \sqrt{b^2 +4}),\\
    \label{beta_identity_1}
    &\beta = \frac{q_{n-1}}{q_n} + (-1)^n\frac{\beta^{n+1}}{q_n},
    \\\label{beta_identity_2}&
    q_n = \frac{1}{\sqrt{b^2+4}}\left(\beta^{-(n+1)} - (-\beta)^{n+1}\right),
    \\\label{beta_identity_3}&
    \sum_{j=0}^{\infty} b\cdot \beta^{2j+1} = 1.
\end{align}

\subsection{Shifted Sudler products}\label{shifted_sudler}
We follow a decomposition approach that was implicitly used already in \cite{grepstad_neum} to prove $\liminf\limits_{N \to \infty}P_N(\phi) > 0$ and more explicitly, in later literature \cite{other_aist_borda,quantum_invariants,tech_zaf}. Here, the Sudler product $P_N(\beta)$ is
decomposed into a finite product with factors of the form

 \begin{equation}
    \label{shifted_sudler_beta}
P_{q_n}(\beta,\varepsilon) := \prod_{r=1}^{q_n} 2 \left\lvert \sin\Big(\pi\Big(r\beta + (-1)^{n}\frac{\varepsilon}{q_n}\Big)\Big)\right\rvert
\end{equation}
where $q_n = q_n(\beta) < N$.
We will first study the case where $\beta = \phi$, which is the easiest case, and will discuss 
the general case afterwards.
Due to the well--known fact that $q_n(\phi) = F_{n+1}$, we will work with shifted Sudler products of the form
\[P_{F_n}(\phi,\varepsilon) := \prod_{r=1}^{F_n} 2 \left\lvert \sin\Big(\pi\Big(r\phi + (-1)^{n-1}\frac{\varepsilon}{F_n}\Big)\Big)\right\rvert.\]
Let $N$ be a positive integer and let $n \gs 2$ be the unique integer such that $F_n \ls N < F_{n+1}$.
The Ostrowski expansion (also called Zeckendorff expansion for the case of the Golden Ratio) simplifies to

\[N = F_n + F_{n_k} + \ldots + F_{n_1}\]
for some appropriate $k,n_1,\ldots,n_k$ where $n_k \ls n-2, n_i - n_{i-1} \gs 2, n_1 \gs 2.$
Denoting\\ \mbox{$N_i = F_n + F_{n_k} + \ldots + F_{n_{i+1}} \text{ for } 1 \ls i \ls k-1$}, we can write

\begin{align}
    P_N(\phi) &= 
    P_{F_n}(\phi)\cdot \prod_{i=1}^k\prod_{r = 1}^{F_{n_i}} 2 \lvert \sin \pi(r\phi + N_i\phi)\rvert
    \\&= P_{F_n}(\phi) \cdot \prod_{i=1}^kP_{F_{n_i}}(\phi,\varepsilon_i),\label{decomp_phi}
\end{align}
where $\varepsilon_i = (-1)^{n_i+1}F_{n_i}N_i\phi$.
Using \eqref{beta_identity_1}, we obtain 
\begin{equation}\label{eps_explicit}\varepsilon_i =F_{n_i}\phi^{n_i}\big((-1)^{n_{i+1}-n_i}\phi^{n_{i+1} - n_i}
+ \ldots (-1)^{n_k -n_i}\phi^{n_k - n_i} + (-1)^{n-n_i}\phi^{n-n_i}\big).\end{equation}
By \eqref{beta_identity_2} and  $n_i \gs 2$, we have \[F_{n_i}\phi^{n_i} = \frac{1}{\sqrt{5}}\big(1 + (-1)^{n_i+1}\phi^{2n_i}\big)
\ls \frac{1 + \phi^6}{\sqrt{5}},\] so we get by 
\eqref{beta_identity_3}
the upper bound

\begin{equation}\label{phi_eps_upper}\varepsilon_i \ls \frac{1}{\sqrt{5}}\left(\left(\phi^2 + \phi^4 + \ldots\right)\right)(1 + \phi^6) = \frac{\phi}{\sqrt{5}}(1+\phi^6) < 0.3 \end{equation}
and the lower bound

\begin{equation}\label{phi_eps_lower}\varepsilon_i \gs -\frac{1}{\sqrt{5}}\left(\left(\phi^3 + \phi^5 + \ldots\right)\right)(1 + \phi^6) = \frac{\phi^2}{\sqrt{5}}(1+\phi^6) > -0.19.\end{equation}
We see that we have reduced the problem of finding a reasonable lower bound on $P_N(\phi)$ to the problem of finding lower bounds on $P_{F_i}(\phi,\varepsilon)$ for $i = 1, \ldots, k$ and $\varepsilon$ in a given interval.\vspace{5mm}

Generalizing this idea to the case of quadratic irrationals with continued fraction period length $1$, we let
\[N = \sum_{i=0}^n b_{i+1}q_i(\beta)\] be the Ostrowski expansion of a positive integer $q_n \ls N < q_{n+1}$.
Following the argument in the proof of \cite[Lemma 3]{tech_zaf}, we 
define 
\[M_{i,a_i} = \sum_{j=n-i+1}^{n} b_{j+1}q_j + a_{n-i}q_i, \quad\quad i = 0,\ldots, n, \quad a_i = 0,\ldots ,b_{n+i-1},\]
to obtain 

\begin{align}
P_N(\beta) & = \prod_{i=0}^{n}\prod_{a_i= 0}^{b_{n-i+1}-1} \prod_{r = M_{i,a_i}}^{M_i,a_i} 2\lvert \sin(\pi r \beta) \rvert
= \prod_{i=0}^{n}\prod_{a_i= 0}^{b_{n-i+1}-1} \prod_{r = 1}^{q_{n-i}} 2\lvert \sin(\pi (M_{i,a_i} + r) \beta) \rvert
\\& = \prod_{i=0}^{n}\prod_{a_i= 0}^{b_{n-i+1}-1} P_{q_{n-i}}(\beta,\varepsilon_{i,a_i})\label{prod_prod_beta}
\end{align}
where

\begin{align}
\label{form_of_eps}\frac{(-1)^{n-i}{\varepsilon_{i,a_i}}}{q_{n-i}} &= M_{i,a_i}\beta,
 \end{align}
 and an empty product is understood to be $1$.
Using the identities \eqref{beta_identity_1} and \eqref{beta_identity_2} leads to

\begin{equation}\label{estim_of_eps}\varepsilon_{i,a_i} = \frac{1}{\sqrt{b^2 + 4}}\left(a_i - b_{n-i+2}\beta + b_{n-i+3}\beta^2 - ...
    + (-1)^ib_{n+1}\beta^{i}\right)(1 + r_i),\end{equation}
where
\begin{equation}\label{size_of_delta}r_{i} = (-1)^{n-i}\beta^{2(n-i+1)}.\end{equation}
As in the Golden Ratio case, we use \eqref{beta_identity_3} to get the bounds 

\begin{equation}\label{eps_upper_bound_quadr}\varepsilon_{i,a_i} \ls \frac{1}{\sqrt{b^2+4}}\left(a_i + b\left(\beta^2 + \beta^4 + \ldots\right)\right)(1 + r_i)\ls \frac{(b-1)+ \beta}{\sqrt{b^2+4}}(1+\beta^2)\end{equation}
and
\begin{equation}\label{eps_lower_bound_quadr}\varepsilon_{i,a_i} \gs \frac{1}{\sqrt{b^2+4}}\left(a_i - b\left(\beta + \beta^3 + \ldots\right)\right)(1 + r_i) \gs -\frac{1}{\sqrt{b^2+4}}(1+\beta^2).\end{equation}
So again we have reduced finding bounds on $P_N(\beta)$ to the problem of finding bounds for
$P_{q_k}(\beta,\varepsilon)$ where $\varepsilon$ lies in a specified interval. The tools to obtain these bounds will be developed in the next section.

\subsection{Limit function}\label{limit_section}
As mentioned in the introduction, Mestel and Verschueren \cite{mestel} showed that 
$P_{F_n}(\phi) = P_{F_n}(\phi,0)$ converges to some positive constant $C_1$, a result 
generalized by Aistleitner, Technau and Zafeiropoulus by the following theorem:\vspace{5mm}

\noindent
{\bf Theorem B} (Aistleitner, Technau and Zafeiropoulus \cite[Theorem 4]{tech_zaf}).
Let $b$ be a positive integer and let $\beta = \beta(b) = [0;\overline{b}\,]$. 
For every $\varepsilon \in \mathbb{R}$, the limit \[G_{\beta}(\varepsilon) := \lim_{n \to \infty}P_{q_n}(\beta,\varepsilon)\]
exists. The convergence is uniform on compact intervals.\vspace{5mm}

The authors also give a rather long, closed expression for $G_{\beta}$ which helped them to approximately compute $G_{\beta}$ and in particular, to calculate $\lim_{N \to \infty} P_{q_n}(\beta) = G_{\beta}(0)$. 

In a recent paper, Aistleitner and Borda \cite{quantum_invariants} quantified this convergence. Their result holds for arbitrary quadratic irrationals $\alpha$, but we state the theorem only for 
irrationals whose continued fraction expansion has period length $1$:\vspace{0.5cm}
\noindent

{\bf Theorem C} (Aistleitner, Borda \cite[Theorem 4]{quantum_invariants}).
Let $\beta(b) := [0;\overline{b}\,]$ %
and let 
\begin{equation}\label{limit_function}
    G_{\beta(b)}(\varepsilon) = 2\pi\lvert \varepsilon + \tfrac{1}{\sqrt{b^2 + 4}}\rvert\prod_{n=1}^{\infty}
     g_n(\beta,\varepsilon) 
\end{equation}
where
\[g_n(\beta,\varepsilon) := \Bigg\lvert\left(1 - \frac{1}{\sqrt{b^2 + 4}}\frac{\left\{n\beta\right\} - \frac{1}{2}}{n}\right)^2 - \frac{\left(\varepsilon + \frac{1}{2\sqrt{b^2 +4}}\right)^2}{n^2}\Bigg\rvert.\] Furthermore,
let $I \subset \mathbb{R}$ be a compact interval. %
Then

\begin{equation}\label{convergence_rate}
    P_{q_k}(\beta,\varepsilon) = \left(1 + \mathcal{O}\left(q_k^{-1/2}\log^{3/4} q_k\right)\right)G_{\beta}(\varepsilon) + \mathcal{O}(q_k^{-2})
\end{equation}
with implied constants only depending on $I$ and $\beta$.\vspace{5mm}\par

With a more detailed analysis, we will show that the multiplicative error in \eqref{convergence_rate} can be decreased further
for arbitrary quadratic irrationals. 

\begin{lem}\label{quantitative_convergence}
Let $\beta, k, G, I$ be as in Theorem C.
 Then we have

\begin{equation}\label{lem_3_estimate}
    P_{q_k}(\beta,\varepsilon) = G_{\beta}(\varepsilon)\left(1 + \mathcal{O}(q_k^{-2/3}\log^{2/3} q_k) + \mathcal{O}(q_k^{-2})\right)
\end{equation}
with implied constants only depending on $I$ and $\beta$ which are explicitly computable.
\end{lem}

We note that our method of proof is not restricted to quadratic irrationals of the form $[0;\overline{\beta}\,]$; taking other quadratic irrationals leads to other implied constants, but the order of convergence remains the same. For this reason, we can improve the convergence rate $\mathcal{O}\left(q_k^{-1/2}\log^{3/4} q_k\right)$ from \eqref{convergence_rate} to
$\mathcal{O}(q_k^{-2/3}\log^{2/3} q_k)$ 
for arbitrary quadratic irrationals.\vspace{2mm}

As a corollary of Theorem C respectively Lemma 4,
we can find for any given $\eta,\gamma > 0$ and any %
interval $I$ an explicitly computable $K_0$ such that
for all $k \gs K_0$ and any $\varepsilon \in I$
\begin{equation}\label{mult_convergence}P_{q_k}(\beta,\varepsilon) > (1 - \eta)G(\beta,\varepsilon) - \gamma .\end{equation}

This implies that

\begin{equation}\label{P_overdef}\overline{P}_{\beta}(\varepsilon) = \min \left\{\min_{0 \ls k \ls K_0} P_{q_k}(\beta,\varepsilon), (1 - \eta)G(\beta,\varepsilon) - \gamma\right\}\end{equation}
fulfills the property $P_{q_k}(\beta,\varepsilon) \gs \overline{P}_{\beta}(\varepsilon)$
for all $k \gs 0$.
However, we can still not take $\overline{P}_{\beta}$ as our computable lower bound due to the following two reasons:
\begin{itemize}
    \item We do not know the exact value of the perturbation $\varepsilon$, but only an interval $[a,b]$ where $\varepsilon$ lies in.
    \item $G$ is defined as an infinite product and is therefore not exactly computable.
\end{itemize}

To solve the first issue, we can show that $G$ and $P_{q_k}$ are both pseudo-concave on zero-free intervals. As a reminder, we call a function $f$ to be pseudo-concave on an interval $[a,b]$ if it fulfills
\begin{equation}\label{pseudo_con}f(x) \gs \min\{f(a),f(b)\}\text{ for any } x \in [a,b].\end{equation}
It is an easy exercise to check that $\log$-concavity implies pseudo-concavity, a fact that was used 
in \cite{tech_zaf} to show that $G$ is pseudo-concave on zero-free intervals.
Similarly, we see by elementary calculus that for any $k \in \mathbb{N}$

\begin{equation}
\pderiv[2]{}{\varepsilon} \log P_{q_k}(\beta,\varepsilon)
= - \sum_{\ell = 0}^{q_k}\frac{\left((-1)^{k}\frac{\pi}{q_k}\right)^2}{\sin^2\left(\pi\big(\ell
\beta + (-1)^k\frac{\varepsilon}{q_k}\big)\right)} < 0
\end{equation}
as long as $\varepsilon$ is no zero of $P_{q_k}$. Hence,
also $P_{q_k}$ is pseudo-concave on zero-free intervals.
Pseudo-concavity is stable under scaling, adding constants and taking a minimum, so we can conclude that \eqref{pseudo_con} also holds for $\overline{P}_{\beta}$.
Therefore, we can estimate

\begin{equation}P_{q_k}(\varepsilon) \gs \overline{P}_{\beta}(\varepsilon) \gs  \min\{\overline{P}_{\beta}(a),\overline{P}_{\beta}(b)\}.\end{equation}

Concerning the second issue, we see that it suffices to compute a finite truncation of the infinite product to obtain a sufficiently small approximation error as shown by the following lemma.

\begin{lem}\label{truncation_lem}
Let $G_{\beta}$ be defined as in \eqref{limit_function}, $T \in \mathbb{N}$,
$I\subseteq \mathbb{R}$ a compact interval.
Writing
\[G_{\beta,T}(\varepsilon) = 2\pi\left\lvert \varepsilon + \tfrac{1}{\sqrt{b^2 + 4}}\right\rvert\prod_{n = 1}^{T} g_n(\beta,\varepsilon),\]
we have for sufficiently large $T$ that
\[G_{\beta}(\varepsilon) = G_{\beta,T}(\varepsilon)\cdot\big( 1 + \mathcal{O}\big(\tfrac{\log T}{T}\big)\big)\]
where the implied constant only depends on $\beta$ and $I$ and is explicitly computable.
\end{lem}

So we find for any $\eta' > 0$ and any zero-free interval $I$ an explicitly computable positive integer $T$ such that for any $\varepsilon \in I$
\begin{equation}\label{mult_convergence_2}%
G_{\beta}(\varepsilon) > (1 -\eta')G_{\beta,T}(\varepsilon).\end{equation}

Defining \begin{equation}\label{Pdef}P_{\beta}(\varepsilon) = \min \{\min_{0 \ls k \ls K_0} P_{q_k}(\beta,\varepsilon), (1-\delta)G_{\beta,T}(\varepsilon) - \gamma\},\end{equation}
where $(1-\delta) = (1 - \eta)(1-\eta')$,
we get from the discussion on pseudo-concavity that
\[P_{q_k}(\beta,\varepsilon) \gs \min\{{P}_{\beta}(a),{P}_{\beta}(b)\}\]
for any $k \gs 1$, with ${P}_{\beta}$ being computable in finitely many steps.
In some cases, this estimate will turn out to be too coarse, since $P_{q_0}(\beta,\varepsilon)$ takes much smaller values than 
$P_{q_n}(\beta,\varepsilon)$ for $n \gs 1$. Hence, we define analogously to \eqref{Pdef} the function

\begin{equation}\label{P*def}P^{*}_{\beta}(\varepsilon) = \min \{\min_{1 \ls k \ls K_0} P_{q_k}(\beta,\varepsilon), (1-\delta)G_{\beta,T}(\varepsilon) - \gamma\}.\end{equation}
Using the same arguments as for ${P}_{\beta}$, we obtain 
\[P_{q_k}(\beta,\varepsilon) \gs \min\{P_{\beta}^{*}(a),P_{\beta}^{*}(b)\}\]
for any $k \gs 1$. 
We will apply these considerations to $\beta(b)$ where $1 \ls b \ls 5$.
Note that ${P}_{\beta},{P}_{\beta}^{*}$ still implicitly depend on several parameters; choosing $\delta,\gamma$
very small leads to rather large values for $K_0,T$ and thus, to larger computational effort to calculate $P,P^{*}$. In the following corollary, we balance this
choice in a way such that $\delta$ and $\gamma$ are small enough to obtain good enough lower bounds
to prove Theorem \ref{main_3} and \ref{main_1}, but are large enough such that $K_0$ and $T$ remain small enough to compute ${P}_{\beta},{P}_{\beta}^{*}$ with reasonable computational effort. To avoid repeating the arguments five times, we will only consider the easiest case, that is the case of the Golden Ratio with $b=1$, for seeing the main structure, as well as the most delicate case, that is, $b=5$.
The intermediate cases where $b = 2,3,4$ can be treated analogously.

\begin{corollary}\label{quantitative_convergence_cor}
Let $G_{\beta,T}$ be defined as in Lemma \ref{truncation_lem}.
Then we have the following.
\begin{enumerate}
    \item[(i)] Let $k \gs 25, \varepsilon \in [a,b] \subseteq [-0.19,0.3]$, $T = 100.000$. Then %
\begin{equation}\label{explicit_phi}P_{F_{k}}(\phi,\varepsilon) \gs 0.94 \cdot \min\{G_{\phi,T}(a),G_{\phi,T}(b)\} - 0.001.\end{equation}
    \item[(ii)] Let $b = 5, k \gs 10, \varepsilon \in [a,b] \subseteq [-0.15,0.93]$, $T = 100.000$. Then
    \begin{equation}\label{explicit_beta5}P_{q_{k}}(\beta,\varepsilon) \gs 0.998 \cdot \min\{G_{\beta,T}(a),G_{\beta,T}(b)\} - 0.0001.\end{equation}
\end{enumerate}
\end{corollary}

Note that the intervals are chosen in a way such that every perturbation $\varepsilon$ possibly appearing in the decompositions
in Section \ref{shifted_sudler}
lies in the corresponding interval, which follows from \eqref{phi_eps_upper} and \eqref{phi_eps_lower} for the Golden Ratio respectively 
\eqref{eps_upper_bound_quadr} and \eqref{eps_lower_bound_quadr} in the case $b = 5$.
To avoid breaking the flow of the paper, we postpone the rather long proofs of Lemma \ref{quantitative_convergence}, Lemma \ref{truncation_lem} and Corollary \ref{quantitative_convergence_cor} to Section \ref{quant_conv_proof}.

\section{Proofs for the case of the Golden Ratio}\label{golden_ratio_sec}
We will prove Theorem \ref{main_3} and Theorem \ref{main_1} together. In fact, we prove the following three statements for $\beta = \beta(b) = [0;\overline{b}\,], 1 \ls b \ls 5$.\\
$\bullet$ For $n \gs 1$, we have
      \begin{equation}\label{ineq1}P_{q_n}(\beta) > P_1(\beta).\end{equation}
$\bullet$
      For $n \gs 1$ and $q_n \ls N < q_{n+1}$, we have
      \begin{equation}\label{ineq2}P_{q_n}(\beta) \ls P_N(\beta).\end{equation}
$\bullet$
    For $n \gs 1$ and $q_n \ls N < q_{n+1}$, we have
    \begin{equation}\label{ineq3}\frac{P_N(\beta)}{N} \ls \frac{P_{q_{n+1}-1}(\beta)}{q_{n+1}-1} .\end{equation}
      
Theorem \ref{main_3} follows immediately from \eqref{limsup_mirror}, \eqref{ineq2} and \eqref{ineq3}
whereas the statement \eqref{thm_1_2} of Theorem \ref{main_1} follows from \eqref{ineq2} and \eqref{ineq3}
and \eqref{min_at1} follows from \eqref{ineq1} and \eqref{ineq2}. Note that in \eqref{ineq3}, we actually prove something stronger than needed in the theorems: for Theorem \ref{main_3}, we only need \eqref{ineq3} to hold for sufficiently large $N$, and for Theorem \ref{main_1} it would suffice to have
$P_N(\beta) \ls P_{q_{n+1}-1}(\beta)$. 

Readers who are familiar with the reflection principle, which is the reason behind the connection between the behaviour of $\liminf_{N \to \infty} P_N(\beta)$ and $\limsup_{N \to \infty} \frac{P_N(\beta)}{N}$)
(see \cite[proof of Lemma 1]{erdos_szekeres} or more explicit in \cite{tech_zaf})
might assume that one can extend \eqref{ineq3} to 
   \begin{equation}\frac{P_{q_{n}}(\beta)}{q_{n}} \ls \frac{P_N(\beta)}{N} \ls \frac{P_{q_{n+1}-1}(\beta)}{q_{n+1}-1} .\end{equation}
This can indeed be proven for the case of the Golden Ratio, but interestingly fails to hold for $b = 5$,
even when assuming $n$ sufficiently large. One can show for example that for and $n\gs 1$ and
$N = q_n + q_{n-1}$, we have $\frac{P_N(\beta)}{N} < \frac{P_{q_{n}}(\beta)}{q_{n}}$.
This can be proven with the tools provided in this paper, but will not be elaborated further.\vspace{5mm}

Returning to the proofs of \eqref{ineq1}\,--\,\eqref{ineq3}, we will concentrate first on the case where $b =1$, that is the case of the Golden Ratio. Since $F_{n+1} = q_n$, we assume for easier notation $F_n < N < F_{n+1}$ and $n \gs 2$. We verify \eqref{ineq1}\,--\,\eqref{ineq3} for all $n \ls 7$ and hence, $N < F_8 = 21$, by explicit computations, so we will assume throughout the rest of this section that $n \gs 8,N \gs 22$.\vspace{5mm}

We start with proving \eqref{ineq1} as a warm-up.
Let $n \gs 3$, then we have 
$P_{F_n}(\phi) \gs P^{*}_{\phi}(0)$
where $P^{*}$ is defined as in \eqref{P*def} with parameters $K_0,T,\delta,\gamma$ determined by Corollary \ref{quantitative_convergence_cor}. During the rest of this paper, we will always mean by $P,P^{*}$
the functions with these fixed parameters, also suppressing the dependence of $P,P^{*}$ on $\phi$ respectively $\beta$ for the sake of readability.
Numerical evaluation yields
$P^{*}(0) \approx 2.22 > 1.86 \approx P_1(\phi,0)$ which proves \eqref{ineq1}.\vspace{4mm}

   \begin{figure}
    \centering
 \includegraphics[width=15cm]{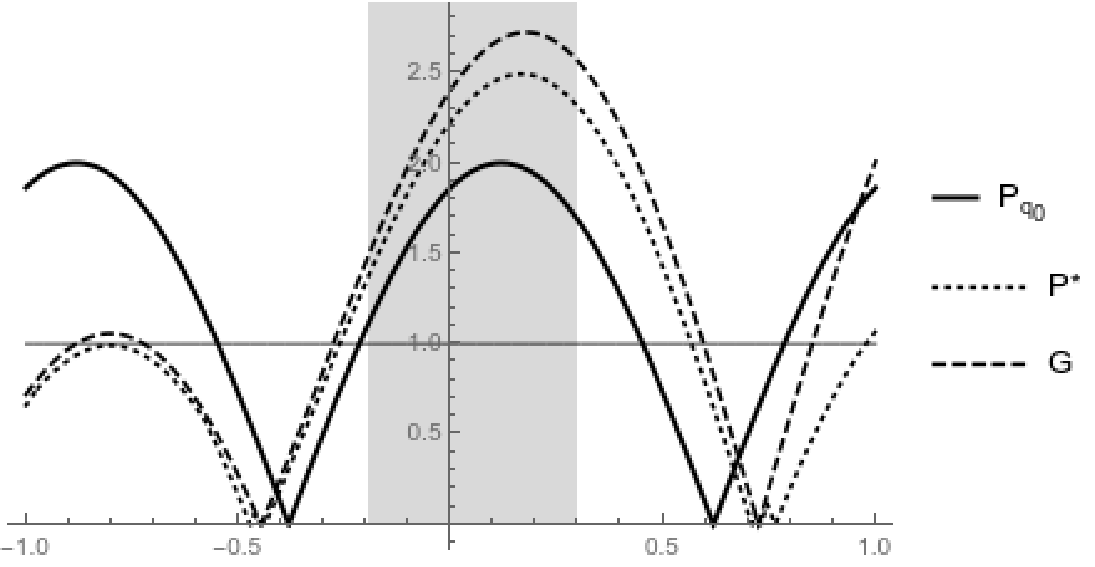}
    \caption{Plot of $P_{q_0},P^{*},G$ for the Golden Ratio case.
    We see that all three functions exceed $1$ in the relevant range $-0.19 < \varepsilon < 0.3$ (grey shading). Looking at $P^{*}$ instead of $P$ gives us a significantly
    better lower bound than $P$ to estimate factors of the form $P_{F_i}(\phi,\varepsilon), i \gs 3$.
    We clearly see that $P^{*}(0) > P_{q_0}(\phi,0)$ which is the reason why \eqref{ineq1} holds.
    }
    \label{fig_golden_ratio}
\end{figure}

To prove \eqref{ineq2}, let $N = F_n + F_{n_k} + \ldots + F_{n_1}$ be the Zeckendorff expansion of $N$.
By the discussion in Section \ref{shifted_sudler} we can write
\begin{align}
\label{left_ineq_prod}
P_N(\phi) =
    &P_{F_n}(\phi) \cdot
    \prod_{i =1}^{k}P_{F_{n_i}}(\phi,\varepsilon_i)
\end{align}
where $-0.19 < \varepsilon_i < 0.3$.
We obtain by actual computation that 
\[\min\{P(-0.19), P(0.3)\} \approx 1.13 > 1,\] so we can deduce that
$P_{F_{n_i}}(\phi,\varepsilon_i) > 1$, and thus, $P_N(\phi) \gs P_{F_n}(\phi)$ follows immediately.\vspace{5mm}

We proceed to show \eqref{ineq3}, which turns out to be more involved.
By applying \eqref{beta_identity_1}, we can write

\begin{align}\label{mirror}P_N(\phi) = \frac{P_{F_{n+1}-1}(\phi)}{\prod\limits_{\ell=N+1}^{F_{n+1}-1}2\lvert\sin(2\pi(\ell\phi))\rvert}
= \frac{P_{F_{n+1}-1}(\phi)}{\prod\limits_{\ell=1}^{F_{n+1}-N-1}2\lvert\sin(\pi(F_{n+1} - \ell)\phi)\rvert}
= \frac{P_{F_{n+1}-1}(\phi)}{\prod\limits_{\ell=1}^{F_{n+1}-N-1}2\lvert\sin(\pi(\ell\phi + (-\phi)^{n+1}))\rvert},
\end{align}
so we see that \eqref{ineq3} is equivalent to showing that 

\begin{equation}
    \label{worst_case_denom}
\prod\limits_{\ell=1}^{F_{n+1}-N-1}2\lvert\sin(\pi(\ell\phi + (-\phi)^{n+1}))\rvert > \frac{F_{n+1}-1}{N}.
\end{equation}
Since \eqref{ineq3} holds trivially for $N = F_{n+1}-1$, we can write

\[F_{n+1} - (N + 1) = F_{n_k} + \ldots + F_{n_1}\]
in its Zeckendorff expansion with $\;2 \ls k \ls n-2$.
Using the arguments from Section \ref{shifted_sudler}, 
with $N$ substituted by $F_{n+1} - (N+1)$ and the additional perturbation coming from $(-\phi)^{n+1}$,
we see that
\begin{align}
    \prod\limits_{\ell=1}^{F_{n+1}-N-1}2\lvert\sin(\pi(\ell\phi + (-\phi)^{n+1}))\rvert =
    \prod_{i=1}^{k} P_{F_{n_i}}(\tilde{\varepsilon}_i)
\end{align}
where 

\begin{equation}\label{perturbed_eps}\tilde{\varepsilon}_i = %
\frac{1}{\sqrt{5}}\big((-1)^{n_{i+1}-n_i}\phi^{n_{i+1} - n_i}
+ \ldots +(-1)^{n_k -n_i}\phi^{n_k - n_i} - (-1)^{n+1-n_i}\phi^{n+1-n_i}\big)(1 + r_i)\end{equation}
with $\lvert r_i \rvert \ls \phi^6$.
The only difference to \eqref{eps_explicit} is the opposite sign in the last summand above. However, it can be shown that we still have $-0.19 < \tilde{\varepsilon}_i< 0.3$ for $i = 1, \ldots, k-1$:
by \eqref{beta_identity_3},
we have
\[\phi^{n+1-n_i} = \phi^{n+1-n_i+1} + \phi^{n+1-n_i+3} + \ldots,\]
so the additional perturbation equals precisely the worst case possible for the estimate on higher orders of $\phi$ appearing in  \eqref{phi_eps_upper} and \eqref{phi_eps_lower}.
Therefore, we can deduce immediately that \mbox{$P_{F_{n_i}}(\phi,\tilde{\varepsilon}_i) > 1$} holds for any $i$.
These arguments are the key for the connection between $\liminf_{N \to \infty} P_N(\phi)$ and $\limsup_{N \to \infty} \frac{P_N(\phi)}{N}$ established in \cite{tech_zaf}. There,
it sufficed to prove that the left-hand side of \eqref{worst_case_denom} is bounded away from $0$, which is not enough to prove \eqref{ineq3}:
in fact, for $n \gs 8$, we have 
\begin{equation}\label{growth_of_fib}\frac{F_{n+1}-1}{N} \ls \frac{F_{n+1}}{F_n} \ls 1.67,\end{equation}
so it suffices to prove

\begin{equation}
    \label{suffices_phi}
\prod_{i=1}^{k} P_{F_{n_i}}(\phi,\tilde{\varepsilon}_i) \gs 1.67.
\end{equation}
 We consider the following case distinction, depending on the Ostrowski coefficients of \mbox{$F_{n+1}-(N+1)$},
 showing that we fulfill \eqref{suffices_phi} in any case.\par

\begin{itemize}
\item Case 1: $k=1$, $n_k=2$.
In this case, we have

\[\prod_{i=1}^{k} P_{F_{n_i}}(\phi,\tilde{\varepsilon}_i) = P_{F_2}(\phi,(-1)^n \phi^{n+1}).\]

Since $\lvert \phi^n \rvert \ls 0.06$ for $n \gs 8$,
we can deduce that
$P_{F_2}(\phi,(-1)^n\phi^n) > 1.69$.\vspace{3mm}

\item Case 2: $k = 1$, $n_k \gs 3$. Observe that for $\tilde{\varepsilon}_k$, the sum in \eqref{perturbed_eps} simplifies to the single term $(-1)^{n-n_k}\phi^{n+1-n_k}$. Thus, we obtain the better estimates

\begin{equation}\label{better_estimate}- 0.12 < - \frac{\phi^3}{\sqrt{5}}(1 + \phi^6)\ls \tilde{\varepsilon_{k}} \ls \frac{\phi^2}{\sqrt{5}}(1 + \phi^6) < 0.19,\end{equation}

which imply
\begin{equation}\label{value_single_l}P_{F_{n_k}}(\phi, \varepsilon_k)\gs \min\{P^{*}(-0.12),P^{*}(0.19)\} > 1.75.\end{equation}

\item Case 3: $k \gs 2$, $n_{k-1} \gs 3$. %
We have for any $n_i \gs 3$ that
\begin{equation}\label{explicit_P*} P_{F_{n_i}}(\phi,\tilde{\varepsilon}_i) \gs \min\{P^{*}(-0.19),P^{*}(0.3)\} \gs 1.395.\end{equation}
Since $1.395\cdot 1.395 > 1.67$ and all other factors $P_{F_{n_i}}(\phi,\tilde{\varepsilon}_i)$ are bounded from below by $1$, we obtain \eqref{suffices_phi}.
\vspace{3mm}

\item
Case 4: $k \gs 2$, $n_{k-1} = 2$. %
Directly bounding $P_{F_{n_2}}(\phi,\tilde{\varepsilon}_2)$ by \mbox{$\min\{P^{*}(-0.19),P^{*}(0.3)\} \approx 1.395$} and $P_{F_{2}}(\phi,\tilde{\varepsilon}_1)$ by $\min\{P(-0.19),P(0.3)\} \approx 1.13$ is not enough to obtain a value that exceeds $1.67$, so we have to argue differently. %
First note that $n_{k-1} = 2$ implies $k =2$.
Since $n \gs 8$, we have $n-n_2 \gs 4$ or $n_2 \gs 5$, so in both cases at least one of the perturbations $\tilde{\varepsilon}_1,\tilde{\varepsilon}_2$ is very close to $0$, improving the estimates sufficiently:

If $n-n_2 \gs 4$, we have
\[-0.07 < - \frac{\phi^4}{\sqrt{5}}\left(1 + \phi^6\right) \ls \tilde{\varepsilon}_2 \ls  \frac{\phi^4}{\sqrt{5}}\left(1 + \phi^6\right) < 0.07\]
and by $\min\{P^{*}(-0.07),P^{*}(0.07)\} > 1.97$, the result follows.

If $n_2 \gs 5$, we see that 
\[-0.15 < -\frac{\phi^3 + \phi^5}{\sqrt{5}}\left(1+ \phi^6\right) \ls \tilde{\varepsilon}_1 \ls \frac{\phi^3 + \phi^5}{\sqrt{5}}\left( 1+ \phi^6\right) < 0.15.\]
Since $P_{F_2}(\phi,\tilde{\varepsilon}_{1}) \gs \min\{P(-0.15),P(0.15)\} > 1.33$,
we get in combination with \eqref{explicit_P*} that 
\[P_{F_2}(\phi,\tilde{\varepsilon}_{1})\cdot P_{F_{n_2}}(\phi,\tilde{\varepsilon}_{2}) > 1.395\cdot 1.33 > 1.67\] which concludes the proof.
\end{itemize}

\section{Proofs for quadratic irrationals of the form $\beta = [0,\overline{b}]$}\label{quadr_irr_sec}

Here we prove \eqref{ineq1}\,--\,\eqref{ineq3} for the case where $b = 5$. The general layout of the proof remains similar to the one for the Golden Ratio, however, we have to be more careful in the numerical analysis. As before, we check numerically with computer assistance that \eqref{ineq1}\,--\,\eqref{ineq3} hold for $n \ls 5$ and $N < q_6 = 18901,$ allowing us to assume from now on that $n \gs 6$.\vspace{4mm}

\begin{figure}
    \centering
 \includegraphics[width=12cm]{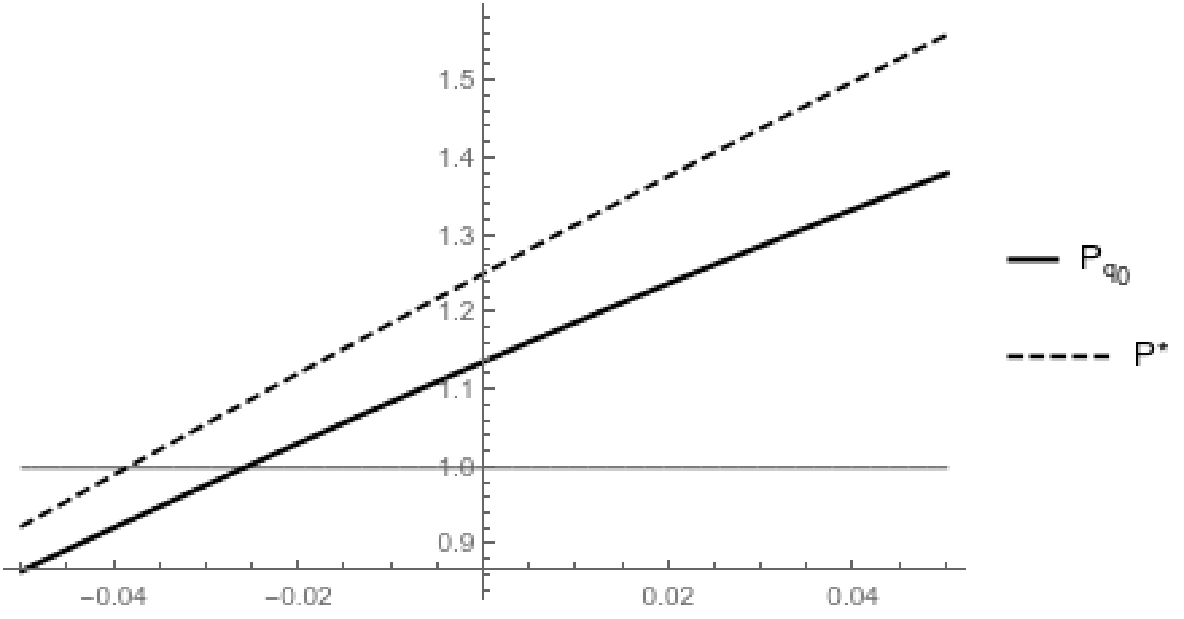}
    \caption{\label{fig_at0} Plot of $P_{q_0},P^{*}$ for the case $b=5$, zoomed in at $\varepsilon = 0$.
    As in the Golden Ratio case, we see that $P_{\beta}^*(0) > P_{q_0}(\beta,0)$, so \eqref{ineq1} also holds for $\beta(5)$.
    }
\end{figure}

Proving \eqref{ineq1} is still straightforward: We have $P_{\beta}^*(0) > P_{1}(\beta,0)$, so \eqref{ineq1}
follows by the same argument as in the case of the Golden Ratio (see Figure \ref{fig_at0}). 

Next, we proceed to show \eqref{ineq2}.
We recall from \eqref{prod_prod_beta} that if
$N = \sum\limits_{i=0}^n b_{i+1}q_i$ is the Ostrowski expansion of $N$,
we can write

\begin{equation}\label{short_decomp}P_N(\beta) = \prod_{i=0}^{n}\prod_{a_i= 0}^{b_{n-i+1}-1} P_{q_{n-i}}(\beta,\varepsilon_{i,a_i})\end{equation}
with 
\begin{equation}\label{estim_of_eps_2}
\varepsilon_{i,a_i} = \frac{1}{\sqrt{b^2 + 4}}\left(a_i - b_{n-i+2}\beta + b_{n-i+3}\beta^2 - ...
    + (-1)^ib_{n+1}\beta^{i}\right)(1 + r_i)\end{equation}
where
$r_{i} = (-1)^{n-i}\beta^{2(n-i+1)}$.
The analysis here is more tedious than in the Golden Ratio case since the limit function does not fulfill $G_{\beta,\varepsilon} > 1$
for all appearing perturbations $\varepsilon$, so in particular, we cannot argue that
 $P^{*}(\varepsilon) > 1$.
       \begin{figure}
    \centering
 \includegraphics[width=15cm]{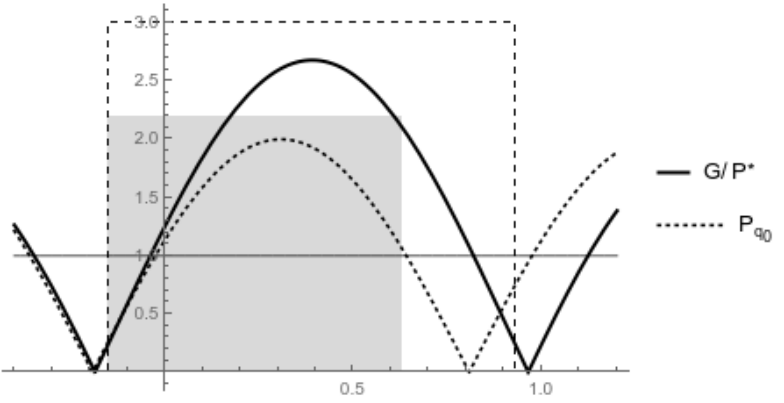}
 \caption{\label{fig_limit_quadr}Plot of $G,P_{q_0},P^{*}$ for the $b=5$ case. Since the convergence 
 of $P_{q_k}$ to $G$ is very fast, even a plot against $P_{q_1}$ or $P^{*}$ is at that scale indistinguishable from the limit function $G$. The grey shaded range indicates the range where a perturbation for $P_1 = P_{q_0}$ needs to be considered, the dashed one for arbitrary $P_{q_k}$.
    We see that $P_{q_k}(\beta,\varepsilon) > 0$ in the considered range, but in contrast to the Golden Ratio case, $P_{q_k}(\beta,\varepsilon)$ can fall significantly below 1, for both $\varepsilon$
    negative and largely positive.
    }
\end{figure}
However, as we see in Figure \ref{fig_limit_quadr}, we only have \mbox{$P^{*}(\varepsilon) < 1$} if $\varepsilon$ is either very small or very large. Looking at \eqref{estim_of_eps_2}, we see that $\varepsilon_{i,a_i}$ is very small if $a_i = 0$ and $b_{n-i+2}$ is ``large'', that is $b_{n-i+2} = 4$ or $b_{n-i+2} = 5$. We will see that
$\varepsilon_{i-1,2},\varepsilon_{i-1,3}$ are values located in an area where
$P^{*}$ is close to its peak, exceeding $1$ by a large margin. As the factors $P_{q_{n-i+1}}(\beta,\varepsilon_{i-1,2}),P_{q_{n-i+1}}(\beta,\varepsilon_{i-1,3})$
are factors of the product 
$\prod\limits_{a_{i-1}= 1}^{b_{n-i+2}-1}P_{q_{n-i+1}}(\beta,\varepsilon_{i-1,a_{i-1}})$ when $b_{n-i+2} = 4$ or $b_{n-i+2} = 5$, these factors will overcompensate the factor
$P_{q_{n-i}}(\beta,\varepsilon_{i,0})$ which might be smaller than $1$,
which results in $\left(\prod\limits_{a_{i-1}= 1}^{b_{n-i+2}-1}P_{q_{n-i+1}}(\beta,\varepsilon_{i-1,a_{i-1}})\right)\cdot
P_{q_{n-i}}(\beta,\varepsilon_{i,0})$ exceeding $1$ in most cases.
Similarly, if the perturbation $\varepsilon_{i,a_i}$ is very large, this implies that $a_i = 4$ and $b_{n-i+1} = 5$, again yielding overcompensating factors that allow us to deduce that the above product exceeds $1$ in most cases. For a detailed explanation and an explicit example of this phenomenon, we refer to \cite[p.28\,--\,30]{tech_zaf}.
Following this observation, we write

\begin{align}\label{orig_prod_prod_3}
P_N(\beta) &= \prod_{i=0}^{n}\prod_{a_i= 0}^{b_{n-i+1}-1} P_{q_{n-i}}(\beta,\varepsilon_{i,a_i}) 
    \\
    &=P_{q_n}(\beta)\cdot 
    \prod_{i=0}^{n}\left(\prod_{a_i = 1}^{b_{n-i+1}-1}
    P_{q_{n-i}}(\beta,\varepsilon_{i,a_i})\right)
    \cdot
    P_{q_{n-(i+1)}}(\beta,\varepsilon_{i+1,0})^{\mathds{1}_{[b_{n-i} \neq 0]}}
\end{align}
with the convention $b_0 = 0$.
For shorter notation, we define
\begin{equation}\label{form_of_Pni}P_{n,i} = \left(\prod_{a_i=1}^{b_{n-i+1}-1} P_{q_{n-i}}(\beta,\varepsilon_{i,a_i})\right)\cdot P_{q_{n-(i+1)}}(\beta,\varepsilon_{i+1,0})^{\mathds{1}_{[b_{n-i} \neq 0]}}.\end{equation}
It follows immediately that showing $P_{q_n}(\beta) \ls P_N(\beta)$ is
equivalent to proving 
\begin{equation}\label{prod_Pni}\prod\limits_{i=0}^n P_{n,i} \gs 1.\end{equation}

    To prove that, we will show that most factors
    fulfill $P_{n,i} \gs 1$. If this cannot be shown for some fixed $i$, we will prove that there
    is another factor $P_{n,j_i}, j_i \neq j_k$ for $i \neq k$ which is much larger than $1$ such that $P_{n,j_i}P_{n,i} \gs 1$.
    This approach is very similar to the proof of \cite[Lemma 3]{tech_zaf}.
    However, in \cite{tech_zaf} it suffices to show that $P_{n,i}>1$ (or $P_{n,i}P_{n,j}>1$) holds whenever $n$ is sufficiently large and for all $i\ls n-N_0$ where $N_0$ is another, arbitrary large constant.  As we cannot make these assumptions here, we will work with the functions $P^{*}$ and $P$ and not with arbitrary close approximations to the limit function $G$.
    Since $P \ls  P^{*} < G$, we need to distinguish even more cases than already considered in the delicate analysis in \cite{tech_zaf}. %
    Note that estimates on $P_{n,i}$ need a special treatment when $n-i$ is small:
    We might not be able to use $P^{*}$ as lower bound, but need to consider the weaker bound $P$.
    Additionally, the error $r_i = q_{n-i}\beta^{n-i+1} - \frac{1}{\sqrt{b^2 +4}}$ appearing in the estimate for $\varepsilon_{i,a_i}$  (see \eqref{estim_of_eps_2}) gets larger the smaller $n-i$ is.
    Therefore, we need to distinguish several cases depending on whether $n-i \gs 2$ or not.\vspace{2mm}
    
   \subsection*{Estimates for $\bf{P_{n.i},\;  n-i \gs 2}$}\par
    In view of \eqref{estim_of_eps_2} and \eqref{form_of_Pni}, we have perturbations of the form
    $\varepsilon_{i,a_i},\varepsilon_{i+1,0}$ for $a_i = 1,\ldots, b_{n-i+1}-1$
    with
\begin{equation}\varepsilon_{i,a_i} = \frac{1}{\sqrt{b^2 + 4}}\left(a_i - b_{n-i+2}\beta + b_{n-i+3}\beta^2 - ...
    + (-1)^ib_{n+1}\beta^{i}\right)(1 + r_i)\end{equation}
where
\begin{equation}r_{i} = (-1)^{n-i}\beta^{2(n-i+1)}.\end{equation}
Since $i \ls n-2$, %
we can deduce that
\begin{equation}\label{size_of_r_i}-\beta^4 \ls r_i \ls \beta^6.\end{equation} 
If $\varepsilon_{i,a_i} > 0$, we will use the estimate $r_i \ls \beta^6$ for upper bounds and
$r_i \gs -\beta^4$ for lower bounds. If $\varepsilon_{i,a_i} < 0$, we use the converse inequalities for upper respectively lower bounds.
Note that these values are very small, so they will not play a main role.
    Using \eqref{beta_identity_3}, %
    we obtain

\begin{equation}\label{coarse_eps_up}a_i - b_{n-i+2}\beta + b_{n-i+3}\beta^2 - ...
    + (-1)^ib_{n+1}\beta^{i}
    \ls a_i + 5\beta^2 + 5 \beta^4 + 5\beta^6 + \ldots  = a_i + \beta.
    \end{equation}\vspace{2mm}
    Since $b_{n-i+1}-1 \gs a_i \gs 1$, we have in particular $b_{n-i+1} \neq 0$, and thus, we know by the properties of the Ostrowski expansion that
    $b_{n-i+2} \ls 4$. This leads to%
    
    \begin{equation}\label{coarse_eps_low}a_i - b_{n-i+2}\beta + b_{n-i+3}\beta^2 - ...
    + (-1)^ib_{n+1}\beta^{i} \gs a_i - 4\beta - 5\beta^3 - 5\beta^5 - \ldots
    = a_i-4\beta - \beta^2.\end{equation}
      Since $P^{*}$ is a lower bound on the pseudo-convex function $P_{q_{n-i}}(\beta,\cdot)$, we can deduce from
     \eqref{size_of_r_i}, \eqref{coarse_eps_up} and \eqref{coarse_eps_low} that if $a_i \gs 1$,
      
        \begin{equation}\label{triv_bound_i}P_{q_{n-i}}(\beta,\varepsilon_{i,a_i})
    \gs \min\Big\{P^{*}\Big(\frac{a_i-4\beta -\beta^2}{\sqrt{29}}(1 - \beta^4)\Big),P^{*}\Big(\frac{a_i+\beta}{\sqrt{29}}(1 + \beta^6)\Big)\Big\}.\end{equation}
    Concerning $\varepsilon_{i+1,0}$, we see by similar arguments that
    
      \begin{equation}\label{triv_bound_i+1}P_{q_{n-(i+1)}}(\beta,\varepsilon_{i+1,0})
    \gs \min\Big\{P^{*}\Big(\frac{-b_{n-i+1}\beta - \beta^2}{\sqrt{29}}(1 + \beta^6)\Big),P^{*}\Big(\frac{- (b_{n-i+1}-1)\beta}{\sqrt{29}}(1 + r_i)\Big)\Big\}\end{equation}
    with $r_i$ being $-\beta^4$ or $\beta^6$, depending on the value of $b_{n-i+1}$.
            Looking at \eqref{form_of_Pni}, we observe that the value of $P_{n,i}$ heavily depends on the coefficients $b_{n-i+1},b_{n-i}$, so we will make a case distinction depending on the values of those coefficients.
    The values on the right-hand sides of \eqref{triv_bound_i} and \eqref{triv_bound_i+1} are used as lower bounds in most cases stated below, so we collect these lower bounds for $P_{q_{n-i}}(\varepsilon_{i,a_i}), a_i = 1,\ldots, 4$ and
    $P_{q_{n-(i+1)}}(\varepsilon_{i+1,0})$ %
    in the following table.
    
        \begin{table}[h]
 \begin{tabular}{||c | c |c||} 
 \hline
 k &If $b_{n-i+1} \neq 0$, then $P_{q_{n-i}}(\beta,\varepsilon_{i,k})\gs$ & If $b_{n-i+1} = k$, then  $P_{q_{n-(i+1)}}(\beta,\varepsilon_{i+1,0})\gs$\\ [0.5ex] 
 \hline\hline
  0 & - &1.20 \\ 
 \hline
 1 &  1.47& 0.97\\
 \hline
 2 & 2.37 & 0.73\\
 \hline
 3 &  2.25 & 0.48\\
 \hline
 4 & 1.12 &  0.24\\
  \hline
\end{tabular}
\vspace{3mm}
\caption{\label{table_(i)} Lower bounds on $P_{q_{n-i}}(\varepsilon_{i,a_i}), P_{q_{n-(i+1)}}(\varepsilon_{i+1,0})$ for $n-i \gs 2$. Note that 
$P_{q_{n-i}}(\varepsilon_{i,0})$ does not appear as a factor in $P_{n,i}$, regardless of the values of 
$b_{n-i+1},b_{n-i}$.}
\end{table}
    \vspace{-2mm}
    
   $\bullet$ Case 1: $b_{n-i} \neq 0, b_{n-i+1} = 0$.
   In this case, $P_{n,i} = P_{q_{n-(i+1)}}(\varepsilon_{i+1,0})$. From Table \ref{table_(i)} we obtain
   $P_{n,i} \gs 1.20$.\vspace{5mm}
   
   $\bullet$ Case 2: $b_{n-i} \neq 0, b_{n-i+1} = 1$.
As in Case 1, we have 
   $P_{n,i} = P_{q_{n-(i+1)}}(\beta,\varepsilon_{i+1,0})$, but the lower bound for $\varepsilon_{i+1,0}$ from \eqref{triv_bound_i+1} is not enough to show
   $P_{n,i} \gs 1$.
   Fixing the values of $b_{n-i+2},b_{n-i+3}$, we obtain the sharper estimate
   \begin{equation}\label{case2_sharper}\varepsilon_{i+1,0} \gs \frac{-\beta + b_{n-i+2}\beta^2 - b_{n-i+3}\beta^3 - \beta^4}{\sqrt{29}}(1 + \beta^6)\end{equation}
   with equality to the bound from \eqref{triv_bound_i+1} in case of $b_{n-i+2} = 0, b_{n-i+3} = 5$.
   We will consider two subcases, depending on $b_{n-i+2},b_{n-i+3}$:
   \begin{itemize}
   \item[--] Case 2a: $b_{n-i+2} \gs 1$ or $b_{n-i+2} = 0, b_{n-i+3} \ls 1$: In this case, \eqref{case2_sharper} is good enough to deduce $P_{n,i}> 1.008$.\vspace{5mm}
   \item[--] Case 2b: $b_{n-i+2} = 0, b_{n-i+3} \gs 2$: 
   Here, \eqref{case2_sharper} does not improve enough over \eqref{triv_bound_i+1} to deduce $P_{n,i} > 1$, so we can only prove that $P_{n,i} > 0.97$.
   We will see after the case distinctions that we can overcompensate
   this factor $P_{n,i}$ with $P_{n,i-2}$. Note that $b_{n-i+3} \neq 0$ implies $i \gs 2$, so 
   $P_{n,i-2}$ is part of the product considered in \eqref{prod_Pni}.
   \end{itemize}\vspace{5mm}
   
   $\bullet$ Case 3: $b_{n-i} \neq 0, b_{n-i+1} = 2$.
   In this case we have $P_{n,i} = P_{q_{n-i}}(\beta,\varepsilon_{i,1})\cdot P_{q_{n-(i+1)}}(\beta,\varepsilon_{i+1,0})$, so by Table \ref{table_(i)} we have
   $P_{n,i} \gs 1.47\cdot 0.73 \gs 1.07.$\vspace{5mm}
   
   $\bullet$ Case 4: $b_{n-i} \neq 0, b_{n-i+1} = 3$.
    This gives $P_{n,i} = P_{q_{n-i}}(\beta,\varepsilon_{i,1})\cdot P_{q_{n-i}}(\beta,\varepsilon_{i,2})\cdot P_{q_{n-(i+1)}}(\beta,\varepsilon_{i+1,0})$, which leads to
    $P_{n,i} \gs 1.47\cdot 2.37\cdot 0.48 \gs 1.67$.\vspace{5mm}

$\bullet$ Case 5: $b_{n-i} \neq 0, b_{n-i+1} = 4$. 
We obtain $P_{n,i} = P_{q_{n-i}}(\beta,\varepsilon_{i,1})\cdot P_{q_{n-i}}(\beta,\varepsilon_{i,2})\cdot P_{q_{n-i}}(\beta,\varepsilon_{i,3})\cdot P_{q_{n-(i+1)}}(\beta,\varepsilon_{i+1,0})$ and so we can deduce
$P_{n,i} \gs 1.47\cdot 2.37\cdot 2.25\cdot 0.24 \gs 1.88$.\vspace{5mm}

$\bullet$ Case 6: $b_{n-i} = 0$. Note that this case covers $b_{n-i+1} = 5$ completely since
 $b_{n-i+1} = 5$ implies $b_{n-1} = 0$.
 \begin{itemize}
 \item[--]Case 6a: $b_{n-i+1} \ls 1$: Then the product for $P_{n,i}$ is empty and hence, $P_{n,i} = 1$.\vspace{5mm}
\item[--]Case 6b: $b_{n-i+1} = 2$. $P_{n,i} = P_{q_{n-i}}(\beta,\varepsilon_{i,1})$, so we obtain $P_{n,1} \gs 1.47$.\vspace{5mm}
 \item[--]Case 6c: $b_{n-i+1} = 3$. $P_{n,i} =  P_{q_{n-i}}(\beta,\varepsilon_{i,1})\cdot P_{q_{n-i}}(\beta,\varepsilon_{i,2})$ and thus, $P_{n,i} \gs 1.47\cdot 2.37 \gs 3.48$.\vspace{5mm}
 \item[--]Case 6d: $b_{n-i+1} = 4$. \mbox{$P_{n,i} = P_{q_{n-i}}(\beta,\varepsilon_{i,1})\cdot P_{q_{n-i}}(\beta,\varepsilon_{i,2})\cdot P_{q_{n-i}}(\beta,\varepsilon_{i,3})$} gives\\
 \mbox{$P_{n,i} \gs 1.47\cdot 2.37 \cdot 2.25\gs 7.83$}.\vspace{5mm}
 \item[--]Case 6e: $b_{n-i+1} = 5$. $P_{n,i} = P_{q_{n-i}}(\beta,\varepsilon_{i,1})\cdot P_{q_{n-i}}(\beta,\varepsilon_{i,2})\cdot P_{q_{n-i}}(\beta,\varepsilon_{i,3})\cdot P_{q_{n-i}}(\beta,\varepsilon_{i,4})$
 which leads to 
  $P_{n,i} \gs  1.47\cdot 2.37 \cdot 2.25 \cdot 1.12 \gs 8.77.$
 \end{itemize}
 \vspace{5mm}
To conclude, we see that the only case where we cannot show that $P_{n,i} > 1$ is Case 2b, which forces the factor $P_{n,i-2}$ to be 
in one of the Cases 6b\,--\,6e. We see that this implies $P_{n,i-2} \gs 1.47$, so we deduce
$P_{n,i}\cdot P_{n,i-2} \gs 1.47\cdot 0.97 > 1$ therefore,
$\prod\limits_{i=0}^{n-2} P_{n,i} \gs 1$ follows.

\subsection*{Estimates for $\bf{P_{n,n-1}}$}\par
For $i = n-1$, we have
\[P_{n,n-1} = \left(\prod_{j=1}^{b_2-1}P_{q_1}(\beta,\varepsilon_{n-1,j})\right)P_{q_0}(\beta,\varepsilon_{n,0})^{\mathds{1}_{[b_1 \neq 0]}} \gs \left(\prod_{j=1}^{b_2-1}P^*(\varepsilon_{n-1,j})\right)P_{q_0}(\beta,\varepsilon_{n,0})^{\mathds{1}_{[b_1 \neq 0]}}.\]
If $b_1 = 0$, we can apply Case $6$ from above and are done.
If, however, $b_1 \neq 0$, we cannot apply the estimate $P_{q_0} \gs P^{*}$, but will compute $P_{q_0}$ directly. We can take the estimates for $\varepsilon_{n,0}$ from \eqref{triv_bound_i+1} with the adjustment that
$r_i$ can now only be bounded from above by $\beta^2$. 
The correspondingly changed values in the table look as follows:

        \begin{table}[h]
 \begin{tabular}{||c |c||} 
 \hline
 k & If $b_2 = k$, then $P_{q_0}(\beta,\varepsilon_{n,0})\gs$\\ [0.5ex]
 \hline\hline
  0 & 1.10 \\ 
 \hline
 1 & 0.89\\
 \hline
 2 &  0.68\\
 \hline
 3 &  0.46\\
 \hline
 4 & 0.23\\
  \hline
\end{tabular}
\vspace{3mm}
\caption{\label{table(i)n-1} Lower bounds on $P_{q_0}(\beta,\varepsilon_{n,0})$. We see that the bounds are lower than in the second column of Table \ref{table_(i)}.}
\end{table}

Going through the calculations in Cases $1,3,4,5$ with these slightly changed values from Table \ref{table(i)n-1} shows that this is still sufficient to prove $P_{n,n-1} > 1$.
Concerning Case 2, which was already the most intricate one before, we have to argue differently. 
We see that without making additional assumptions, we can only prove that $P_{n,n-1} > 0.89$.
As above, we will
find a factor $P_{n,j}$ which overcompensates the factor $P_{n,n-1}$ and differs from the factors used to overcompensate Case 2b for any $i \ls n-2$. However, we cannot fix the $j_i$ as in the case $i \ls n-2$ with $j_i =i-2$, but $j_i$ will depend on possibly all Ostrowski coefficients. 
We consider the following cases:\vspace{5mm}
\begin{itemize}
\item[--]Case $2a'$: There exists some $k \in \{3,\ldots, n-1\}$ such that $b_k = 0$. Taking the smallest $k$ with that property, we know that $b_{k-1} \neq 0$. This leads for $2 \ls j_{n-1} := n-k+1 \ls n-2$ to Case 1 where we have proven that
$P_{n,j_{n-1}} \gs 1.20$.
Thus, \mbox{$P_{n,j_{n-1}}\cdot P_{n,n-1} \gs 1.20 \cdot 0.89 > 1.06$}.\vspace{5mm}

\item[--]Case $2b'$: $b_k \neq 0$ for all $k \in \{3,\ldots, n-1\}$ and %
there exists some $k \in \{3,\ldots, n-1\}$ such that $b_k \gs 3$.
This implies that for $2 \ls j_{n-1} := n-k+1 \ls n-2$ we are either in Case 4 or Case 5 where we have shown that
$P_{n,j_{n-1}} \gs 1.67$, so again we can prove $P_{n,j_{n-1}}\cdot P_{n,n-1} \gs 1.$\vspace{5mm}

\item[--] Case $2c'$: $b_k \in \{1,2\}$ for all $k \in \{3,\ldots, n-1\}$.
This implies in particular $b_3 \gs 1, b_4 \ls 2$, so we can improve the estimate \eqref{triv_bound_i+1} to
\begin{equation}\label{case2_n-1_improved}\frac{-\beta + \beta^2 - 2\beta^3 -\beta^4}{\sqrt{29}}(1+\beta^2) < \varepsilon_{n,0} < 0.\end{equation}
This leads to the sharper estimate $P_{n,n-1} = P_{q_0}(\beta,\varepsilon_{n,0}) > 0.96$. 
By the assumption $b_k \in \{1,2\}$, we see that $P_{n,n-2}$ is in Case 2 or Case 3. 
If $P_{n,n-2}$ is in Case 3, we obtain from the discussion above that $P_{n,n-2} > 1.07$ which suffices 
to show $P_{n,n-2}\cdot P_{n,n-1} > 1.$ If we are in Case 2, we use the improved bound \eqref{case2_n-1_improved} (with $1+\beta^6$ instead of $1+\beta^2$) to get
$P_{q_{2}}(\beta,\varepsilon_{n-1,0}) > 1.045$, again deducing $P_{n,n-2}\cdot P_{n,n-1} > 1$.
\end{itemize}\vspace{4mm}

As mentioned, the factor $P_{n,j_{n-1}}$ appearing in Cases $2a',2b',2c'$ is not used by some $i\neq n-1$ to overcompensate Case 2b:
Overcompensating factors for Case 2b always come from Case 6, whereas $P_{n,j_{n-1}}$ is always in one of the Cases 1\,--\,5.

\subsection*{Estimates for $\bf{P_{n,n}}$}\par
Observe that by definition, $b_0 = 0, b_1 < 5$, so here we are always in Case 6a\,--\,6d.
Case 6a is an empty product, which is defined to be equal to $1$, so we can concentrate on the case where $2 \ls b_1 \ls 4$.
Using estimates as in \eqref{triv_bound_i} and \eqref{triv_bound_i+1}, we obtain

\[0 \ls \varepsilon_{n,k} \ls (3\beta+\beta^2)(1 + \beta^2), \quad k \in \{1,2,3\}\]
which implies that
$P_1(\beta,\varepsilon_{n,k}) \gs 1.13$.
Therefore, we have $P_{n,n} = \prod\limits_{k=1}^{b_1-1} P_1(\beta,\varepsilon_{n,k}) \gs 1$.

Summing up all estimates on $P_{n,i}$ for $0 \ls i \ls n$, we obtain \eqref{prod_Pni}, and thus, \mbox{$P_{q_n}(\beta) \ls P_N(\beta)$} follows.\phantom{\qedhere}

\begin{proof}[Proof of \eqref{ineq3}]

We proceed to prove that for all $q_n \ls N < q_{n+1}$ we have

\begin{equation}\label{(ii)}\frac{P_N(\beta)}{N} \ls \frac{P_{q_{n+1}-1}(\beta)}{q_{n+1}-1}.\end{equation}
To do so, we apply the same procedure as done in \eqref{worst_case_denom} for the Golden Ratio to see that \eqref{ineq3} is equivalent to

\begin{equation}\label{thm3_suffice}
    \prod_{\ell = 1}^{q_{n+1} - N - 1}2\lvert \sin(\pi (\ell\beta + (-\beta)^{n+2})\rvert
    \gs \frac{q_{n+1}-1}{N}.
\end{equation}

Let
$\sum\limits_{i = 0}^{n} b_{i+1}q_i$ be the Ostrowski expansion of $q_{n+1}- N - 1$, possibly with 
leading zero coefficients.
By an analogous argument as in \eqref{orig_prod_prod_3}, we obtain

    \begin{equation}\label{prod_prod_3}
    \prod_{\ell = 1}^{q_{n+1} - N - 1}2\lvert \sin(\pi (\ell\beta + (-\beta)^{n+2})\rvert=
    P_{q_n}(\tilde{\varepsilon}_{0,0})^{\mathds{1}_{[b_{n+1} \neq 0]}}
    \cdot
    \prod_{i=0}^{n}\prod_{a_i = 1}^{b_{n-i+1}-1} P_{q_{n-i}}(\tilde{\varepsilon}_{i,a_i})
    \cdot P_{q_{n-(i+1)}}(\tilde{\varepsilon}_{i+1,0})^{\mathds{1}_{[b_{n-i} \neq 0]}}
    \end{equation}
where

\begin{align}
\label{form_of_eps_tilde}\frac{(-1)^{n-i}{\tilde{\varepsilon}_{i,a_i}}}{q_{n-i}} &= M_{i,a_i}\beta + (-\beta)^{n+2}.
 \end{align}
Accordingly, we get

\begin{equation}\label{estim_of_eps_tilde}\tilde{\varepsilon}_{i,a_i} = \frac{1}{\sqrt{b^2 + 4}}\left(a_i - b_{n-i+2}\beta + b_{n-i+3}\beta^2 - ...
    + (-1)^ib_{n+1}\beta^{i} - (-1)^{i+1}\beta^{i+1}\right)(1 + r_i),\end{equation}
which shows us that we can interpret $\tilde{\varepsilon}_{i,a_i}$ as the ``normal'' $\varepsilon_{i,a_i}$ defined in \eqref{estim_of_eps} with an additional negative Ostrowski coefficient $b_{n+2} = -1$. 
This additional perturbation will only contribute significantly when $i$ is very small. 
This leads to the following lemma:

\begin{lem}\label{from2_tilde_lem}
Let
\[\tilde{P}_{n,i} := \prod_{a_i = 1}^{b_{n-i+1}-1} P_{q_{n-i}}(\tilde{\varepsilon}_{i,a_i})
    \cdot P_{q_{n-(i+1)}}(\tilde{\varepsilon}_{i+1,0})^{\mathds{1}_{[b_{n-i} \neq 0]}},\quad i = 0,\ldots, n\]
    with $n,b_i,\tilde{\varepsilon}$ as before and assume that $b_{n+1} \neq 0$ or $b_n \neq 0$.
    Then we have
\begin{equation}
    \label{prod_from2}
\prod_{i=2}^n \tilde{P}_{n,i} \gs 0.97 \text { and } \tilde{P}_{n,1} \gs 1.
\end{equation}
\end{lem}

\begin{proof}
    If $i \gs 1$, we see analogously to the Golden Ratio case that the additional perturbation vanishes in the estimate for the geometric series: by \eqref{beta_identity_3}, we have
$\beta^{i+1} = 5\beta^{i+2} + 5 \beta^{i+4} + \ldots,$
so the additional perturbation equals precisely the worst-case higher-order estimates used in  \eqref{coarse_eps_up} and \eqref{coarse_eps_low}. Therefore,  
we can use Table \ref{table_(i)} for $5 \ls i \ls n-2$ to deduce by the same arguments as in \eqref{ineq2} that we have $\tilde{P}_{n,i} > 1$ or, if Case 2b applies, $\tilde{P}_{n,i}\cdot \tilde{P}_{n,i-2} > 1$.
Note that if $i \in \{3,4\}$, we cannot consider $\tilde{P}_{n,i-2}$ to show \eqref{prod_from2}, hence 
we simply take the estimate $P_{n,i} \gs 0.97$ if Case 2b applies. We observe that by the assumptions made in Case 2b, at most one of $\tilde{P}_{n,3},\tilde{P}_{n,4}$ can be in this case. This factor is bounded from below by $0.97$, with the other factor exceeding $1$. Since the special cases $i \in \{1,n-1,n\}$
can be treated completely analogously to the proof of \eqref{ineq2}, we conclude the proof of the lemma.\phantom{\qedhere}
\end{proof}
    
    Returning to the proof of \eqref{ineq3}, we argue similarly to the case of the Golden Ratio: we need to obtain factors being significantly larger than $1$
    to gain the extra factor for $\frac{q_{n+1}-1}{N}$ in
    \eqref{thm3_suffice}. 
Bounding $\frac{q_{n+1}-1}{N}$ by \[\frac{q_{n+1}-1}{N} \ls \frac{q_{n+1}}{q_n} \approx 5+ \beta\]
turns out to be too coarse, so we will use a better bound, depending on the Ostrowski coefficients of $q_{n+1}-(N+1)$.

\begin{proposition}\label{better_estim_thm3}
Let $n \gs 2$ and $\sum\limits_{i = 0}^{n} b_{i+1}q_i$ be the Ostrowski expansion of $q_{n+1}-(N+1)$, possibly with 
leading zero coefficients.
Then
\begin{equation}\label{nontrivial_thm3}
    \frac{q_{n+1}-1}{N} \ls \frac{27}{25 - (5b_{n+1} + b_n)} =: u(b_n,b_{n+1}).
\end{equation}
\end{proposition}

\begin{proof}
By the properties of the Ostrowski expansion, we obtain
\[q_{n+1} - (N+1) < b_{n+1}q_n + (b_n +1)q_{n-1},\]
which implies by the recursion formula from \eqref{recursions}
that

\[N \gs (5 - b_{n+1})q_n + b_nq_{n-1} \gs \left(25 - (5b_{n+1} +b_n)\right)q_{n-1}.\]
Applying \eqref{recursions} again, we get $q_{n+1} \ls 27q_{n-1}$,
and the result follows.\phantom{\qedhere}\vspace{5mm}
\end{proof}

    Leaving the possibility of $b_{n+1} = b_{n} = 0$ %
    for the moment aside, we will show that the additional factor exceeding $1$ comes from
    $P_{q_n}(\beta,\tilde{\varepsilon}_{0,0})^{\mathds{1}_{[b_{n+1} \neq 0]}}\cdot \tilde{P}_{n,0}\cdot \tilde{P}_{n,1}$:
    we recall from Cases 3\,--\,6 above that a large $b_{n-i+1}$ leads to a large factor $P_{n,i}$, significantly exceeding $1$, a fact which will hold for $\tilde{P}_{n,i}$ as well. So if
    $u(b_n,b_{n+1})$
    is large, we know that $b_{n+1}$ or $b_n$ has to be large, and thus, $\tilde{P}_{n,0}$
    exceeds $1$ significantly. If $u(b_n,b_{n+1})$ is small, then it suffices that $\tilde{P}_{n,0}\cdot \tilde{P}_{n,1}$
    stays just above $1$. %
    To be more specific, we show the following lemma:
    
    \begin{lem}\label{factor_Pn0Pn1}
    Let $n,N,b_i,\tilde{\varepsilon}$ as above and assume we have $b_{n+1} \neq 0$ or
     $b_{n} \neq 0$.
    Then
    \begin{equation}\label{thm3_suffice_1}P_{q_n}(\beta,\tilde{\varepsilon}_{0,0})^{\mathds{1}_{[b_{n+1} \neq 0]}}\cdot\tilde{P}_{n,0}\cdot\tilde{P}_{n,1} \gs u(b_n,b_{n+1}).\end{equation}    
    \end{lem}

\begin{proof}%

Setting $b_{n+2} = -1$, we obtain from \eqref{estim_of_eps}
that for $n \gs 2$,
\[\frac{\beta}{\sqrt{29}}(1 - \beta^4) \ls \tilde{\varepsilon}_{0,0} \ls \frac{\beta}{\sqrt{29}}(1 + \beta^6),\]
which leads to 
\begin{equation}\label{P_0,0}P_{q_n}(\beta,\tilde{\varepsilon}_{0,0}) \gs 1.47.\end{equation}
 Note that
\[\tilde{P}_{n,0} = \prod_{k=1}^{b_{n+1} -1}P_{q_{n}}(\beta,\tilde{\varepsilon}_{0,k})P_{q_{n-1}}(\beta,\tilde{\varepsilon}_{1,0})^{\mathds{1}_{[b_n \neq 0]}}\]
    where 
    \begin{align}\frac{k + \beta}{\sqrt{29}}\left(1 - \beta^{4}\right) &\ls \tilde{\varepsilon}_{0,k} \ls \frac{k + \beta}{\sqrt{29}}\left(1 + \beta^{6}\right),\\
   \frac{-b_{n+1}\beta - \beta^2}{\sqrt{29}}(1 - \beta^4) &\ls  \tilde{\varepsilon}_{1,0} \ls  \frac{-b_{n+1}\beta - \beta^2}{\sqrt{29}}(1 + \beta^6)\label{latter}\end{align}
   with reverted inequality symbols for \eqref{latter} in case of $b_{n+1} = 0$.
Since $N > q_n$, we can deduce that $b_{n+1} \ls 4$ and $b_n = 0$ if $b_{n+1} = 4$.
    Computing lower bounds for $\tilde{P}_{n,0}$ as in Cases 1\,--\,6 above leads to the following table:

    \begin{table}[h]
 \begin{tabular}{||c | c |c | c |c|c|c||} 
 \hline
  k &$P^*(\tilde{\varepsilon}_{0,k})\gs$ & If $b_{n+1} = k$, & If $b_{n+1} = k, b_{n} \neq 0$, & If $b_{n+1} = k, b_{n} = 0$& $u(1,k)\ls$& $u(5,k) \ls$\\
 & &$P^*(\tilde{\varepsilon}_{1,0})\gs$&
 $1.47^{\mathds{1}_{[b_{n+1} \neq 0]}}\cdot \tilde{P}_{n,0}\gs$ & $1.47^{\mathds{1}_{[b_{n+1} \neq 0]}}\cdot \tilde{P}_{n,0}\gs$& &  \\[0.5ex] 
 \hline\hline
 0 & - & 1.20 & 1.20 & -& 1.125 & 1.35\\
 \hline
 1 &  2.37& 0.97 & 1.425 & 1.47 & 1.422 & 1.8\\
 \hline
 2 & 2.66 & 0.73 & 2.54 & 3.48& 1.93 & 2.7\\
 \hline
 3 &  2.25 & 0.48 & 4.44 & 9.26 & 3 & 5.4\\
 \hline
 4 & - &  0.24 & - & 20.85& 6.75 & - \\
  \hline
\end{tabular}\vspace{3mm}
\caption{\label{tab_Pn0} Lower bounds on $P_{q_n}(\beta,\tilde{\varepsilon}_{0,0})^{\mathds{1}_{[b_{n+1} \neq 0]}}\cdot \tilde{P}_{n,0}$. The values in Rows 3 and 4 exceed $u(1,k)$ in any case.}
\end{table}
\vspace{-5mm}

We see that $P_{q_n}(\tilde{\varepsilon}_{0,0})\cdot \tilde{P}_{n,0}\gs u(b_n,b_{n+1})$
holds whenever $b_n \ls 1$, %
but might fail for $2 \ls b_n \ls 5$.
In these cases, we will see that $\tilde{P}_{n,1}$ gives us another factor that exceeds $1$ significantly, closing the gap between $P_{q_n}(\tilde{\varepsilon}_{0,0})^{\mathds{1}_{[b_{n+1} \neq 0]}}\cdot\tilde{P}_{n,0}$ and $u(5,k)$:
  note that $b_n \gs 2$ implies $b_{n+1} \ls 3$, so we obtain that\[\tilde{P}_{n,1} = \left(\prod_{k=1}^{b_n -1}P_{q_{n-1}}(\beta,\tilde{\varepsilon}_{1,k})\right)\cdot P_{q_{n-2}}(\beta,\tilde{\varepsilon}_{2,0})^{\mathds{1}_{[b_{n-1} \neq 0]}}\]
    where  \begin{align}\label{estim_eps_tilde_1}\frac{k - 3\beta - \beta^2}{\sqrt{29}}(1 - \beta^4) \ls & \tilde{\varepsilon}_{1,k} \ls \frac{k - \beta^2}{\sqrt{29}}(1 + \beta^6),\\\label{estim_eps_tilde_2}
    \frac{-b_n\beta + \beta^3}{\sqrt{29}}(1 + \beta^6) \ls & \tilde{\varepsilon}_{2,0} \ls \frac{-b_n\beta + 3\beta^2 + \beta^3}{\sqrt{29}}(1 - \beta^4).
    \end{align}

       \begin{table}
 \begin{tabular}{||c | c |c |c | c||} 
 \hline 
 k &$P^*(\tilde{\varepsilon}_{1,k})\gs$ & If $b_{n} = k$,  & If $b_{n} = k, b_{n-1} \neq 0$ & If $b_{n} = k, b_{n-1} = 0$\\
 & & then $P^*(\tilde{\varepsilon}_{2,0})\gs$&then $\tilde{P}_{n,1}\gs$ &then $\tilde{P}_{n,1}\gs$ \\[0.5ex] 
 \hline\hline
 1 &  1.68 & - & - & -\\
 \hline
 2 & 2.48 & 0.78& 1.31 & 1.68\\
 \hline
 3 &  2.40 & 0.54 & 2.24 &4.16\\
 \hline
 4 & 1.52&  0.29 & 2.89 & 9.99\\
  \hline
  5 & - &  - & - &15.19\\
  \hline
\end{tabular}\vspace{3mm}
\caption{\label{tab_Pn1} Lower bounds on $\tilde{P}_{n,1}$ for $b_n \gs 2$, using the estimates \eqref{estim_eps_tilde_1} and \eqref{estim_eps_tilde_2}. In any case, we have
$\tilde{P}_{n,1} \gs 1.31$.
}
\end{table}

We see from Table \ref{tab_Pn1} that $b_n \gs 2$ implies $\tilde{P}_{n,1} \gs 1.31$.
Combining this with the values obtained in Column 3 of Table \ref{tab_Pn0} and comparing it with $u(5,k)$, 
we conclude the proof of Lemma \ref{factor_Pn0Pn1}.\phantom{\qedhere}
\end{proof}

\subsection*{The special case $\bf{b_n = b_{n+1} = 0}$}
We are left to consider the special case $b_{n+1} = b_n = 0$ which implies $P_{q_n}(\beta,\tilde{\varepsilon}_{0,0})^{\mathds{1}_{[b_{n+1} \neq 0]}}\cdot \tilde{P}_{n,0} =1$.
We set $k = \min\{i \gs 1 : b_{n-i} \neq 0\}$ and assume first that $k \ls n-2$.
This implies that $\tilde{P}_{n,k}$ is in Case 1 and thus,
$\tilde{P}_{n,k} = P_{q_n-(k+1)}(\beta,\varepsilon_{{k+1},0})$.
Since $b_{n-k+1} = b_{n-k+2} = \ldots = b_{n+1} = 0,\; b_{n+2} = -1$, we see that

\[\frac{-\beta^4}{\sqrt{29}}(1+\beta^6) \ls \varepsilon_{{k+1},0} \ls \frac{\beta^3}{\sqrt{29}}(1+\beta^6),\]
a range where we can prove that
$\tilde{P}_{n,k} \gs 1.24$.
Similarly to Lemma \ref{from2_tilde_lem}, we show that

\[\prod_{i=k+1}^n \tilde{P}_{n,i} \gs 0.89\]
with the following adjustments:\vspace{5mm}

$\bullet$ For $i = n-1$, we cannot use the case distinction in Cases $2a',2b',2c'$, so we can only deduce \mbox{$\tilde{P}_{n,n-1} > 0.89$}.\vspace{5mm}

$\bullet$ Recall that in Lemma \ref{from2_tilde_lem} (which treated the case $b_{n+1} \neq 0$ or $b_n \neq 0$), we found no overcompensating factor if $i \in \{3,4\}$.
This translates now to finding no compensating factor $P_{n,i-2}$ if $i \in \{k+1,k+2\}$. However, by construction, $b_{n-k} \neq 0, b_{n-k+1} = 0$ rules out an appearance of Case 2b there, so we do not have to consider an additional factor below $1$.
Since \[\left(\prod_{i=0}^{k-1} \tilde{P}_{n,i}\right)\cdot P_{q_n}(\beta,\tilde{\varepsilon}_{0,0})^{\mathds{1}_{[b_{n+1} \neq 0]}}\]
is an empty product and $1.26 \cdot 0.89 \gs u(0,0) = 1.08$, we can complete the case $k \ls n-2$ and 
are left with the possibility that $k = n-1$. We see that
$\tilde{P}_{n,n-1} = P_{q_0}(\beta,\varepsilon_{n,0})$ with 
\[\frac{-\beta^4}{\sqrt{29}}(1+\beta^2) \ls \varepsilon_{{n},0} \ls \frac{\beta^3}{\sqrt{29}}(1+\beta^2).\]
Using $P_{q_0}$ instead of $P^{*}$, we can compute that
$\tilde{P}_{n,n-1} > 1.13$.
Since the only other non-trivial factor from \eqref{prod_prod_3} is $\tilde{P}_{n,n} > 1$, 
we have finished the last case and therefore, also the proof of \eqref{ineq3}.

\phantom{\qedhere}
\end{proof}

\section{Proofs of Lemma \ref{quantitative_convergence}, Lemma \ref{truncation_lem} and Corollary \ref{quantitative_convergence_cor}}\label{quant_conv_proof}

\subsection*{Proof of Lemma \ref{quantitative_convergence}}
We will follow the main structure of the proof in \cite[Theorem 4]{quantum_invariants}, improving at some points the estimates by refined arguments and making the constants hidden in the \mbox{$\mathcal{O}\text{--notation}$} explicit. Note that we do not aim to determine the best constants possible, but are contented with sufficiently small constants to deduce Corollary \ref{quantitative_convergence_cor} for $\delta,K_0$ small enough to compute $P$ and $P^{*}$ with acceptable computational help.
The following proposition will be used several times in the proof of Lemma \ref{quantitative_convergence} to find reasonable bounds on products whose factors are close to 1.

\begin{lem}\label{exp_log_lemma}
\item
\begin{enumerate}
    \item[(i)] Let $\lvert a_n \rvert \ls \frac{1}{2}$ and $\lvert a_n\rvert \ls \frac{C}{n}, N,M \in \mathbb{N}$. Then 
\begin{equation}
    \label{exp_log_estim}
\prod_{n = N}^{M} (1 + a_n) \gs 1 - \Bigg(\Big\lvert\sum_{n=N}^{M} a_n\Big\rvert + \frac{C^2}{N-1}\Bigg).
\end{equation}
\item[(ii)] If $\lvert a_n \rvert \ls \frac{1}{2}$ and $\Big\lvert\sum\limits_{n=N}^{M} a_n\Big\rvert \ls c \ls \frac{1}2$, then 
\begin{equation}
    \label{exp_log_estim_2}
\prod_{n = N}^{M} (1 + a_n) \ls 1 + (1+c)\Big\lvert\sum_{n=N}^{M} a_n\Big\rvert.
\end{equation}
\item[(iii)] If $0 < a_n < \frac{1}{2}$, then 
\begin{equation}
    \label{exp_log_estim_3}
\prod_{n = N}^{M} (1 - a_n) \gs 1 - 2\sum_{n=N}^{M}a_n.\end{equation}
\end{enumerate}
\end{lem}
\begin{proof}\item
\begin{enumerate}
    \item[(i)]
We make use of the Taylor series expansion
\[\log(1+x) = -\sum_{k =1}^{\infty} \frac{(-x)^k}{k}\]
to show that
\begin{equation}\label{lower_log_bound}\log(1 + a_n) \gs  a_n - \frac{1}{2}\sum_{k=2}^{\infty}\lvert a_n\rvert ^k = a_n - \frac{1}{2}a_n^2\left(\frac{1}{1-\lvert a_n\rvert}\right)
\gs a_n - a_n^2 \gs a_n - \frac{C^2}{n^2}.\end{equation}
Summing over \eqref{lower_log_bound}, we obtain

\begin{equation*}%
    \sum_{n=N}^{M} \log(1 + a_n) \gs \sum_{n=N}^{M} a_n - \frac{C^2}{N-1}.
\end{equation*}
With $e^x \gs 1+x$, we can deduce the desired result.

\item[(ii)]
Using the estimate $\log(1+x) \ls x$ and that for
$0 \ls x \ls c \ls \frac{1}{2}$
\[\exp(x) \ls 1 + x + x^2 \ls 1 + (1+c)x, \]
we obtain
\eqref{exp_log_estim_2}.%

\item[(iii)] Here we use that for $0 < x < \frac{1}{2}$ we have
$\log(1-x) \gs -2x$ and once more, $e^x \gs 1 + x$ to deduce \eqref{exp_log_estim_3}.
\end{enumerate}\phantom{\qedhere}
\end{proof}

\subsection*{Proof overview}

We will prove the inequality 
\begin{equation*}
    P_{q_k}(\beta,\varepsilon) \gs G_{\beta}(\varepsilon)\left(1 + \mathcal{O}\Big(\big(\tfrac{\log q_k}{q_k}\big)^{2/3}\Big)\right) + \mathcal{O}(q_k^{-2})
\end{equation*}
with explicitly specified constants. This is the inequality we need
to prove the main results of this paper.
The other inequality can be shown by very similar methods (with different explicit values of the implied constants). For better readability, we define $f(x) = \lvert 2\sin(\pi x)\rvert$ and $\delta_k = \|q_k\beta\|$. 
Note that $f$ is $1$--periodic and an even function, facts that will be used implicitly several times below.
Taking the last factor of $P_{q_k}(\beta,\varepsilon)$ out of the product and using 
\eqref{alpha_deltak}, %
we obtain 

\begin{equation}\label{last_factor_out}P_{q_k}(\beta,\varepsilon) = f\Big(\delta_k + \frac{\varepsilon}{q_k}\Big)
\prod_{n=1}^{q_k-1} f\Bigg(\frac{np_k}{q_k} + (-1)^k\Big(\Big\{\frac{n}{q_k}\Big\} - \frac{1}{2}\Big) + (-1)^k\frac{2\varepsilon + q_k\delta_k}{2q_k}\Bigg).
\end{equation}
Observe that the factors on the right-hand side of \eqref{last_factor_out} only depend on the residue class of \mbox{$n\!\mod q_k$}.
Using properties \eqref{conv_ident} and \eqref{bijection},
we deduce

\[P_{q_k}(\beta,\varepsilon) = f\Big(\delta_k + \frac{\varepsilon}{q_k}\Big)
\prod_{n=1}^{q_k-1} f\Bigg(\frac{n}{q_k}  - \Big(\Big\{\frac{nq_{k-1}}{q_k}\Big\} - \frac{1}{2}\Big)  - \frac{2\varepsilon + q_k\delta_k}{2q_k}\Bigg).\]
We see that \[\prod_{n =1}^{q_k-1}f\Big(\frac{n}{q_k}\Big) = \prod_{n =1}^{q_k-1}\big\lvert 1 - e^{2\pi i n/q_k}\big\rvert
= \lim_{x \to 1}\frac{x^{q_k}-1}{x -1} = q_k,\]
so we can write

\begin{equation}\label{step_before_identity}P_{q_k}(\beta,\varepsilon) = f\Big(\delta_k + \frac{\varepsilon}{q_k}\Big)\,q_k
\prod_{n=1}^{q_k-1} \frac{f\Big(\frac{n}{q_k}  - \Big(\big\{\frac{nq_{k-1}}{q_k}\big\} - \frac{1}{2}\Big)  - \frac{2\varepsilon + q_k\delta_k}{2q_k}\Big)}{f\big(\frac{n}{q_k}\big)}.\end{equation}
Combining the $n$--th and ($q_k-n$)\,--th factor of \eqref{step_before_identity} for $1 \ls n < \frac{q_k}{2}$ and using the trigonometric identity $f(u-v)f(u+v) = \lvert f^2(u) - f^2(v)\rvert$,
we obtain

\begin{equation}\label{half_product}
    P_{q_k}(\beta,\varepsilon) =
    f(\delta_k + \tfrac{\varepsilon}{q_k})\,q_k
    \prod_{ 0 < n < q_k/2} \left\lvert\frac{f^2\Big(\frac{n}{q_k} - \big(\{\frac{nq_{k-1}}{q_k}\} - \frac{1}{2}\big)\delta_k\Big) - f^2\big(\frac{2\varepsilon + q_k\delta_k}{2q_k}\big)}{f^2\big(\frac{n}{q_k}\big)}\right\rvert
\end{equation}
with an additional factor 
\begin{equation}\label{even_factor}\frac{f\left(\frac{1}{2} - \frac{2\varepsilon + q_k\delta_k}{2q_k}\right)}{f(\frac{1}{2})} = \frac{f\left(\frac{1}{2} - \frac{2\varepsilon + q_k\delta_k}{2q_k}\right)}{2}\end{equation}
if $q_k$ is even.\vspace{5mm}

Next, we consider a monotonically non-decreasing integer-valued function $\psi(t)$ fulfilling \mbox{$\psi(t) = \Omega(\log t)$}, but $\psi(t) = o(t^{4/5})$, which will be optimized later. Furthermore, let $N_0$ be a sufficiently large integer, depending only on $I$ and $\beta$.
 We denote
\[g_n(\varepsilon) = g_n(\beta,\varepsilon) = \Bigg\lvert\Bigg(1 - \frac{1}{\sqrt{b^2 + 4}}\frac{\left\{n\beta\right\} - \frac{1}{2}}{n}\Bigg)^2 - \frac{\left(\varepsilon + \frac{1}{2\sqrt{b^2 +4}}\right)^2}{n^2}\Bigg\rvert\]
and 
\[p_n(k,\varepsilon) = p_n(\beta,k,\varepsilon) = \left\lvert\frac{f^2\Big(\frac{n}{q_k} - \big(\{\frac{nq_{k-1}}{q_k}\} - \frac{1}{2}\big)\delta_k\Big) - f^2\big(\frac{2\varepsilon + q_k\delta_k}{2q_k}\big)}{f^2\big(\frac{n}{q_k}\big)}\right\rvert.\]
We split the factors on the right-hand side of \eqref{half_product} into several groups and prove Lemma \ref{quantitative_convergence} via the following steps:

\begin{enumerate}
\item[$\bullet$ Step 1:]
\begin{equation}\label{subresult_1}\frac{f\left(\frac{1}{2} - \frac{2\varepsilon + q_k\delta_k}{2q_k}\right)}{2} \gs 1 + \mathcal{O}(q_k^{-2}).\end{equation}
    \item[$\bullet$ Step 2:]
    \begin{equation}\label{subresult_2}\prod_{\psi(q_k) < n < \frac{q_k}{2}} p_n(k,\varepsilon) \gs 1 + \mathcal{O}\Big(\frac{\log q_k}{\psi(q_k)}\Big).\end{equation}
     \item[$\bullet$ Step 3:]
     \begin{equation}\label{subresult_3}
     \prod_{n = \psi(q_k)+1}^{\infty} g_n(\varepsilon) \ls 1 + \mathcal{O}\Big(\frac{\log \psi(q_k)}{\psi(q_k)}\Big).\end{equation}
     
    \item[$\bullet$ Step 4:]%
    \begin{equation}\label{subresult_4}\prod\limits_{n=N_0+1}^{\psi(q_k)} p_n(k,\varepsilon) \gs
    \Big(1 + \mathcal{O}\Big(\frac{\psi(q_k)^2}{q_k^2}\Big)\Big){\prod\limits_{n=N_0+1}^{\psi(q_k)} g_n(\varepsilon)}.\end{equation}
    \item[$\bullet$ Step 5:]    
    \begin{equation}\label{subresult_5}\prod\limits_{n=1}^{N_0} p_n(k,\varepsilon) \gs
    \big(1 + \mathcal{O}\big(q_k^{-2}\big)\big){\prod\limits_{n=1}^{N_0} g_n(\varepsilon)} + \mathcal{O}\big(q_k^{-2}\big).\end{equation}
    \item[$\bullet$ Step 6:] %
    \begin{equation}\label{subresult_6}f\Big(\delta_k + \frac{\varepsilon}{q_k}\Big) \gs  2\pi\Big\lvert \varepsilon + \frac{1}{\sqrt{b^2+4}}\Big\rvert\Big(1 + \mathcal{O}\big(q_k^{-2}\big)\Big)
    + \mathcal{O}(q_k^{-2}).\end{equation}
\end{enumerate}

Combining Steps 1\,--\,6%
, we deduce

\begin{equation}\label{lem3_aim}P_{q_k}(\varepsilon) \gs G(\varepsilon)\Big(1 + \mathcal{O}\Big(\frac{\log q_k}{\psi(q_k)}+ \frac{\psi(q_k)^2}{q_k^2}\Big)\Big) + \mathcal{O}\big(q_k^{-2}\big).
\end{equation}

We minimize the error term in \eqref{lem3_aim} by setting $\psi(t) = t^{2/3}\log^{1/3} t$ which leads to the desired result.
The only step where the order of magnitude differs from \cite[Theorem 4]{quantum_invariants} is in Step 4. %
This leads to a different optimal choice of $\psi$ and to a faster convergence rate.

\subsection*{Proof of Step 1}
For better readability, we define 

\[
h(\varepsilon) = \frac{2\varepsilon + q_k\delta_k}{2q_k}, \quad 
c_1 := \max_{k \gs 1}q_k\cdot \max_{\varepsilon \in I}\lvert h(\varepsilon)\rvert.\]
If $k$ is sufficiently large, then
    $\frac{c_1}{q_k} < \frac{1}{2},$
and thus we can estimate 

\[
\frac{f\left(\frac{1}{2} - h(\varepsilon)\right)}{2} \gs 1 - \frac{\lvert h(\varepsilon)\rvert}{2}f'\left(\tfrac{1}{2} + \lvert h(\varepsilon)\rvert\right)
\gs 1 - \frac{\pi^2 c_1^2}{q_k^2} = 1 + \mathcal{O}(q_k^{-2}).
 \]

\subsection*{Proof of Step 2}

If $k$ is large enough such that $\psi(q_k) > c_1+1$, we can remove for $n > \psi(q_k)$ the absolute values of $p_n(k,\varepsilon)$.
We use the trigonometric identity %

\begin{equation}\label{trig_square}\frac{\sin^2(a-b)}{\sin^2(a)}%
= 1 - \sin(2b)\cot(a) + \left(\cot^2(a)-1\right)\sin^2(b), \quad a \neq 0\end{equation}

to get

\begin{align}\label{trig_identity}\nonumber
    &\frac{f^2\Big(\frac{n}{q_k} - \big(\big\{\frac{nq_{k-1}}{q_k}\big\}-\frac{1}{2}\big)\delta_k\Big)}{f^2(n/q_k)}
    =\\& 1 - \sin\Big(2 \pi \big(\big\{\tfrac{nq_{k-1}}{q_k}\big\}-\tfrac{1}{2}\big)\delta_k\Big)\cot\big(\pi\tfrac{n}{q_k}\big)+ \big(\cot^2(\pi\tfrac{n}{q_k}) -1\big)\sin^2\Big(\pi \big(\big\{\tfrac{nq_{k-1}}{q_k}\big\}-\tfrac{1}{2}\big)\delta_k\Big).
\end{align}

Since quadratic irrationals have bounded partial quotients, we see from \eqref{qkqkalpha} that there exists some fixed $c_2(\beta)>0$ such that $\delta_k \ls \frac{c_2}{q_k}$. Therefore, we have
\[2 \pi \big(\big\{\tfrac{nq_{k-1}}{q_k}\big\}-\tfrac{1}{2}\big)\delta_k \ls \frac{\pi c_2}{2q_k},\]
which leads with $\sin^2(x) \ls x^2$ to 

\begin{equation}
    \label{large_p_estim_3}\Big(\cot^2\big(\pi\tfrac{n}{q_k}\big) -1\Big)\sin^2\Big(\pi \big(\big\{\tfrac{nq_{k-1}}{q_k}\big\}-\tfrac{1}{2}\big)\delta_k\Big)
\gs -\Big(\pi \big(\big\{\tfrac{nq_{k-1}}{q_k}\big\}-\tfrac{1}{2}\big)\delta_k\Big)^2 \gs -\frac{c_2^2\pi^2}{4q_k^2}.
\end{equation}
By $x/2 \ls \sin(x) \ls x$ for $x \in \left[0,\tfrac{\pi}{2}\right],$
we get that
\begin{equation}\label{large_p_estim_1}
\frac{f^2\big(\frac{c_1}{q_k}\big)}{f^2\big(\frac{n}{q_k}\big)} \ls \frac{4c_1^2}{n^2}.
\end{equation}
Combining these estimates, we obtain 

\begin{equation*}%
    p_n(k,\varepsilon) \gs 1 - 
    r(n) - \frac{c_3}{n^2}
\end{equation*}
where \[c_3 = c_2\pi^2/4 + 4c_1^2, \quad 
r(n) = \sin\big(2 \pi \big(\big\{\tfrac{nq_{k-1}}{q_k}\big\}-\tfrac{1}{2}\big)\delta_k\big)\cot\big(\pi\tfrac{n}{q_k}\big) \ls \frac{2c_1}{n}.\]
For $k$ large enough, we have $1 - r(n) \gs \frac{1}{2}$,
so

\begin{align}\nonumber
\prod_{\psi(q_k) < n < \frac{q_k}{2}} p_n(k,\varepsilon) &\gs
\prod_{\psi(q_k) < n < \frac{q_k}{2}} \big(1 - r(n) - \frac{c_3}{n^2}\big)
\\&\gs \prod_{\psi(q_k) < n < \frac{q_k}{2}} \big(1 - r(n)\big) \prod_{\psi(q_k) < n < \frac{q_k}{2}}\left(1 - \frac{2c_3}{n^2}\right).\label{r_n_ineq}
\end{align}
By Proposition \ref{exp_log_lemma} (iii), we get

\begin{equation}\label{large_p_estim}
    \prod_{\psi(q_k) < n < \frac{q_k}{2}}\left(1 - \frac{2c_3}{n^2}\right)
    \gs 
    1 - \frac{4c_3}{\psi(q_k)-1}.
\end{equation}
To bound the first product in \eqref{r_n_ineq}, we apply partial summation and employ the monotonicity of $\cot(x)$ on $0 < x < \frac{\pi}{2}$ to obtain
\begin{align}\label{koksma_estimate}\Big\lvert\sum_{\psi(q_k) < n < \frac{q_k}{2}} r(n) \Big\rvert
&\ls \cot(\pi \psi(q_k)/q_k)\max_{\psi(q_k) \ls n < q_k/2} \Bigg\lvert \sum_{n=\psi(q_k)+1}^{q_k} h(x_n) - \int_{0}^{1} h(x)\, \mathrm{d}x \Bigg\rvert
\end{align}
where 
$h(x) = \sin(2\pi(x - \tfrac{1}{2})\delta_k)$ and $x_n = \{n\beta\}$.
A short elementary computation shows that we can bound the variation by
\[V(h) = 2\sin(\pi\delta_k) < \frac{2\pi c_2}{q_k}.\]
Using a standard discrepancy estimate in the form of \cite[Theorem 1.1]{dst}, we have
for $N < q_k$ 
\[N\cdot D_N^{*}(\{\beta\},\{2\beta\},\ldots,\{N\beta\}\} \ls (b+2)\left\lceil \frac{k+3}{2} \right\rceil,\]
where $D_N^{*}$ denotes the star-discrepancy.
Applying Koksma's inequality 
yields

\begin{align}\label{large_p_estim_4}
\Big\lvert\sum_{\psi(q_k) < n < \frac{q_k}{2}} r(n)  \Big\rvert
\ls \frac{(b+2)(k+4)2c_2}{\psi(q_k)}.
\end{align}
For an introduction on discrepancy and Koksma's inequality, we refer the reader to \cite{kuipers}.
Using $r(n) \ls \frac{2c_1}{n}$ and Proposition \ref{exp_log_lemma} (i), we obtain
\begin{equation}\label{large_p_estim_5}\prod_{\psi(q_k) < n < \frac{q_k}{2}} \big(1 - r(n)\big) 
\gs 1 - \Big(\frac{(b+2)(k+4)2c_2}{\psi(q_k)} + \frac{4c_1^2}{\psi(q_k) -1}\Big).
\end{equation}
Combining \eqref{large_p_estim} and \eqref{large_p_estim_5}, we get the lower bound

\begin{equation}\label{large_p_final}
    \prod_{\psi(q_k) < n < \frac{q_k}{2}} p_n(k,\varepsilon) 
    \gs 1 - \frac{(b+2)(k+4)c_2 + 4c_1^2 + 4c_3}{\psi(q_k)-1} = 1 + \mathcal{O}\Big(\frac{\log q_k}{\psi(q_k)}\Big).
\end{equation}

\subsection*{Proof of Step 3}

Observe that for $k$ sufficiently large and $n \gs \psi(q_k)$, one can remove the absolute values of 
\[g_n(\varepsilon) = \Bigg\lvert \Bigg(1 - \frac{1}{\sqrt{b^2 + 4}}\left(\frac{\{n\beta\} - \frac{1}{2}}{n}\right)\Bigg)^2
 - \frac{\big(\varepsilon + \frac{1}{2\sqrt{b^2+4}}\big)^2}{n^2}\Bigg\rvert\]
and therefore, 

\[\prod_{n = \psi(q_k)+1}^{\infty}g_n(\varepsilon) \ls \Bigg(\prod_{n = \psi(q_k)+1}^{\infty}\big(1 - \ell(n)\big)\Bigg)^2\]
where
\begin{equation}\label{ell_n}\ell(n) = \frac{1}{\sqrt{b^2 + 4}}\left(\frac{\{n\beta\} - \frac{1}{2}}{n}\right).\end{equation}
Due to a result in Ostrowski's famous work \cite{ostrowski}, we get for
$N = \sum\limits_{i=0}^{\ell} b_{i} q_{i}$ %
the estimate
\begin{equation*}%
\Big\lvert \sum_{n=1}^{N} \{n\beta\} - 1/2 \Big\rvert \ls \frac{b(k+1)}{2}.\end{equation*}
If $b>1$, we use the triangle inequality and $\log(N) \gs k\log(b)$ to obtain %
\[\Big\lvert \sum_{n=M}^{N} \{n\beta\} - 1/2 \Big\rvert %
\ls b\Big(\frac{\log(N)}{\log(b)}+1\Big).\]
Applying summation by parts and $N \to \infty$, we get

\begin{equation}
    \label{sum_k_bound}
    \Big\lvert\sum_{n = \psi(q_k)+1}^{\infty} \ell(n)\Big\rvert \ls \frac{b}{\sqrt{b^2 +4}\log(b)}\frac{\log \psi(q_k)}{\psi(q_k)} + \frac{1}{\psi(q_k)}. 
\end{equation}
For $k$ sufficiently large, the left-hand side of \eqref{sum_k_bound} is bounded by $1/4$, and thus, \mbox{$\text{Proposition \ref{exp_log_lemma} (ii)}$} gives

\begin{align}\nonumber
\prod_{n = \psi(q_k)}^{\infty}g_n(\beta,\varepsilon)
&\ls
\Bigg(1 + \frac{5}{4}\cdot\Big\lvert \sum_{n = \psi(q_k)}^{\infty} \ell(n)\Big\rvert\Bigg)^2
\ls 1 + 3\cdot\Big\lvert \sum_{n = \psi(q_k)}^{\infty} \ell(n)\Big\rvert
\\&\ls
1 + \Bigg(\frac{3b}{\sqrt{b^2 +4}\log(b)}\frac{\log(\psi(q_k) +1)}{\psi(q_k)} + \frac{3}{\psi(q_k)}\Bigg)
= 1 + \mathcal{O}\Big(\frac{\log \psi(q_k)}{\psi(q_k)}\Big).
\label{large_for_g}
\end{align}

For $b =1$ and $k \gs 10$ we use the estimate $q_k \gs \big(\frac{3}{2}\big)^k$. Following the same steps as in the $b>1$ case, we end up with an error term as in \eqref{large_for_g}, but with
$\log(b)$ substituted by $\log(3/2)$. Note that this argument also works for arbitrary quadratic irrationals since Ostrowski's estimates can be generalized to all irrationals with bounded partial quotients.

\subsection*{Proof of Step 4}
Before we start with the actual proof of Step 4, we need a technical, but nevertheless important lemma to
find good bounds for terms of the form $\mathlarger{\frac{\sin^2(x+y) - \sin^2(z)}{\sin^2(x)}}$ under certain conditions on $x,y,z$. This lemma
gives us a power saving in comparison to applying the estimate $\sin(t^2) = t^2(1 + \mathcal{O}(t^2))$ for all terms independently as done in \cite[Theorem 4]{quantum_invariants}, leading to the improvement of the convergence rate in Lemma \ref{quantitative_convergence}.%

\begin{lem}\label{better_sin}
Let $0 <  \max\{\lvert y\rvert,\lvert z\rvert\} \ls x  < \frac{1}{2}$
and
\begin{equation}\label{requirement_1}2z^2 \ls (x+y)^2,\end{equation}
\begin{equation}\label{requirement_2}\frac{(x+y)^2- z^2}{x^2} \gs \frac{1}{2}.\end{equation}
Then we have
\begin{equation}\label{sin_approx}
  \frac{\sin^2(x+y) - \sin^2(z)}{\sin^2(x)} 
\gs \frac{(x+y)^2 - z^2}{x^2}\Bigg(1 - \Big(2x\lvert y \rvert + \frac{4}{3}y^2 +  \frac{6x^6}{7!} + x^5\lvert y\rvert +z^2\Big)\Bigg).
\end{equation}

\end{lem}

\begin{proof}
By Taylor approximates and $x \gs \lvert y \rvert$, we have
\begin{align*}
\sin(x) &\ls x\Big(1 - \frac{x^2}{3!} + \frac{x^4}{5!}\Big),\\%
\sin(x+y) &\gs (x+y)\Big(1 - \frac{(x+y)^2}{3!} + \frac{(x+y)^4}{5!}
- \frac{(x+y)^6}{7!}\Big) \\&\gs (x+y)\Big(1 - \frac{x^2}{3!} + \frac{x^4}{5!}\Big)
- (x+y)\Big(\frac{x\lvert y \rvert}{3} + \frac{y^2}{6} + \frac{x^6}{7!} + \frac{2^6 x^5\lvert y\rvert}{7!}\Big).
\end{align*}
Since $x \ls \frac{1}{2}$
implies that $\frac{1}{1-\frac{x^2}{6} + \frac{x^4}{5!}} < \frac{3}{2}$,
we obtain

\begin{equation}\label{tech_lem_1}\frac{\sin(x+y)}{\sin(x)} \gs \frac{x+y}{x}\Bigg(1 - \Big(\frac{x\lvert y\rvert}{2} + \frac{y^2}{3} + \frac{3x^6}{2\cdot 7!} + \frac{x^5\lvert y\rvert}{4}\Big)\Bigg).\end{equation}
We square \eqref{tech_lem_1} and use \eqref{requirement_1} to deduce that

\begin{equation}\label{sin_first_total}
    \frac{\sin^2(x+y)}{\sin^2(x)} %
\gs \frac{(x+y)^2}{x^2}
- \frac{(x+y)^2-z^2}{x^2}
\Big(2x\lvert y\rvert + \frac{4}{3}y^2 + \frac{6x^6}{7!} + x^5\lvert y\rvert\Big).
\end{equation}
Concerning $-\frac{\sin^2(z)}{\sin^2(x)}$, we use Taylor approximations and $x < \frac{1}{2}$ to obtain
\[\sin(x) \gs x\left(1 - \frac{x^2}{6}\right) \gs x\left(\frac{1}{1 + \frac{x^2}{4}}\right)\] 
which leads by $\sin^2(z) \ls z^2$ to %

\begin{equation}\label{second_bound_sin}
    \frac{\sin^2(z)}{\sin^2(x)} \ls \frac{z^2}{x^2}\left(1 + \frac{x^2}{2}\right).
\end{equation}
Note that \eqref{requirement_2} implies

\begin{align}
    \frac{(x+y)^2}{x^2} - \frac{z^2}{x^2}\Big(1 + \frac{x^2}{2}\Big)
    &= \frac{(x+y)^2-z^2}{x^2} -\frac{z^2}{2}
    \gs \frac{(x+y)^2-z^2}{x^2} - \Big(\frac{(x+y)^2}{x^2} - \frac{z^2}{x^2}\Big)z^2
    \\&= \frac{(x+y)^2-z^2}{x^2}\left(1 - z^2\right).\label{third_bound_sin}
\end{align}
Combining \eqref{sin_first_total}, \eqref{second_bound_sin} and \eqref{third_bound_sin}, we get the desired result.\phantom{\qedhere}

\end{proof}

Returning to the proof of Step 4, we set
\[x = x(n) = \pi\frac{n}{q_k},\quad
y = y(n) = -\pi \left(\left\{\tfrac{nq_{k-1}}{q_k}\right\}-\tfrac{1}{2}\right)\delta_k,\quad
z = z(n) =  \pi\frac{2\varepsilon+q_k\delta_k}{2q_k}.\]
Note that $\lvert y \rvert \ls \frac{\pi c_2}{2 q_k}, \; \lvert z \rvert \ls \frac{\pi c_1}{q_k}$
where $c_1,c_2$ are defined as in the proofs of Step 1 and 2.
This implies that for $n \gs N_0$ with $N_0$ chosen sufficiently large, we 
can apply Lemma \ref{better_sin} to obtain

\begin{align}\label{middle_sized_1}\nonumber&\prod_{N_0 < n \ls \psi(q_k)}\frac{f^2\Big(\frac{n}{q_k} - \big(\big\{\frac{nq_{k-1}}{q_k}\big\}-\frac{1}{2}\big)\delta_k\Big)}{f^2(\frac{n}{q_k})}
 - \frac{f^2\left( \frac{2\varepsilon+q_k\delta_k}{2q_k}\right)}{f^2(\frac{n}{q_k})}
\gs \\&\prod_{N_0 < n \ls \psi(q_k)}\frac{(x+y)^2 - z^2}{x^2}\prod_{N_0 < n \ls \psi(q_k)} \left(1 - \Big(c_2\pi^2\frac{n}{q_k^2} + \frac{c_2^2\pi^2}{3}\frac{1}{q_k^2} + \frac{6\pi^6n^6}{7!q_k^6}
+ \frac{\pi^6n^5}{q_k^6}
+ \frac{\pi^2c_1^2}{q_k^2}\Big)\right).
\end{align}
Since %
$n \ls \psi(q_k) = o(q_k^{4/5})$, we have for $k$ sufficiently large that
\[\frac{6\pi^6n^6}{7!q_k^6} + \frac{\pi^6n^5}{q_k^6} \ls \frac{2n}{3q_k^2}.\]
This implies that the error in the last product of \eqref{middle_sized_1} can be bounded by
\begin{equation}\label{step_3_bound_1}\pi^2\frac{(c_2 + 2/3)n + \frac{c_2^2}{3} + c_1^2}{q_k^2}
= \frac{un + v}{q_k^2}
\end{equation}
with $u = \pi^2(c_2 + \frac{2}{3}), v = \pi^2(\frac{c_1^2}{3} + c_1^2)$.
If $k$ is large enough, then the right-hand side of \eqref{step_3_bound_1} can be bounded from above by $1/2$, so  we get by Proposition \ref{exp_log_lemma} (iii) that

\begin{equation}\label{middle_sized_2}\prod_{N_0 < n \ls \psi(q_k)}
\Big( 1 -\frac{un + v}{q_k^2}\Big)
\gs 1 - \Big(u\frac{(\psi(q_k)+1)^2}{q_k^2} + \frac{2v}{3}\frac{\psi(q_k)}{q_k^2}\Big).
\end{equation}
We are left to compare

\[\prod\limits_{N_0 < n < \psi(q_k)}\frac{\lvert (x(n)+y(n))^2 - z(n)^2\rvert}{x(n)^2}
= \prod\limits_{N_0 < n < \psi(q_k)} \Bigg(1 - q_k\delta_k\Big(\Big\{\tfrac{nq_{k-1}}{q_k}\Big\}-\tfrac{1}{2}\Big)\Bigg)^2 - \frac{\big(\varepsilon+\frac{q_k\delta_k}{2}\big)^2}{n^2}
\]
with the product 
\[\prod\limits_{N_0 < n < \psi(q_k)}g_n(\varepsilon) = \prod\limits_{N_0 < n < \psi(q_k)}\left(1 - \frac{1}{\sqrt{b^2 + 4}}\frac{\left\{n\beta\right\} - \frac{1}{2}}{n}\right)^2 - \frac{\left(\varepsilon + \frac{1}{2\sqrt{b^2 +4}}\right)^2}{n^2}.\]
Writing \[s(n) =  1 - \frac{1}{\sqrt{b^2 + 4}}\Bigg(\frac{\{n\beta\} - \frac{1}{2}}{n}\Bigg), 
\quad s'(n) = 1 - q_k\delta_k\left(\left\{\frac{nq_{k-1}}{q_k}\right\}-\frac{1}{2}\right),\]
we see that
\begin{equation}\label{s(n)}s'(n) = s(n) + \underbrace{\Big(\Big\{\frac{nq_{k-1}}{q_k}\Big\}-\frac{1}{2}\Big)\Big( q_k\delta_k - \frac{1}{\sqrt{b^2 + 4}}\Big)}_{=:A} + \underbrace{\Big(\Big\{\frac{nq_{k-1}}{q_k}\Big\}- \left\{n\beta\right\}\Big)\frac{1}{\sqrt{b^2+4}}}_{=:B}.
\end{equation}
From the identities \eqref{beta_identity_1} and \eqref{beta_identity_2}, we can deduce that for $k$ sufficiently large%

\begin{equation}\label{step_3(1)}\lvert A \rvert \ls \frac{\beta^{2(k+1)}}{2\sqrt{b^2+4}}
\ls \frac{1}{(b^2+4)^{3/2}}\frac{1}{q_k^2}
.\end{equation}
To treat $B$, we will show that for sufficiently large $k$, 

\begin{equation}
\Big\lvert\Big\{\frac{nq_{k-1}}{q_k}\Big\}- \big\{n\beta\big\}\Big\rvert
= \left\{n\Big\lvert\frac{q_{k-1}}{q_k}-\beta\Big\rvert\right\},
\end{equation}
that is, the function $\{x\}$ has no discontinuity between $\frac{nq_{k-1}}{q_k}$ and $n\beta$.
Such a ``jump'' would imply that there exists an integer $j$ such that

\begin{equation}\label{contradiction_step3}n\frac{q_{k-1}}{q_k} \gs j \gs n\beta \quad
\text{or} \quad n\beta \gs j \gs n\frac{q_{k-1}}{q_k}.\end{equation}
We assume the former to be the case, with the other one working analogously. Let
\mbox{$q_{\ell} \ls n < q_{\ell+1}$} with $\ell < k-1$. 
By \eqref{qkqkalpha} and \eqref{best_approx},
we obtain

\[n\Big\lVert \frac{nq_{k-1}}{q_k} \Big\rVert \gs q_{\ell}\Big\lVert \frac{q_{\ell+1}q_{k-1}}{q_k} \Big\rVert \gs \frac{q_{\ell}}{(b+2)q_{\ell+1}} \gs \frac{1}{(b+1)(b+2)}.\]
So we have by \eqref{contradiction_step3}

\[\frac{q_{k-1}}{q_k} - \beta \gs\frac{q_{k-1}}{q_k} - \frac{j}{n} \gs
\frac{\lVert nq_{k-1}/q_k \rVert}{n} \gs \frac{1}{(b+1)(b+2)n^2}.
\]
Since $q_{k-1} = p_k$ we know from \eqref{approx_quality} that
\[\frac{q_{k-1}}{q_{k}} - \beta\ls\frac{1}{q_k^2},\]
a contradiction to $n \ls \psi(q_k) = o(q_k^{4/5})$
for $k$ sufficiently large.
 Hence, we can deduce

\begin{equation}\label{no_jump_end}\frac{\{n\beta\}- \{n\frac{q_{k-1}}{q_k}\}}{n} 
\ls \frac{\big\{n\big\lvert\frac{q_{k-1}}{q_k}-\beta\big\rvert\big\}}{n} \ls
\frac{1}{q_k^2}.\end{equation}
Using \eqref{no_jump_end}, combined with \eqref{s(n)} and \eqref{step_3(1)}, we obtain

\begin{equation}\label{s_final}s'(n) \gs s(n) - \left(\frac{1}{(b^2+4)^{3/2}}+ \frac{1}{\sqrt{b^2+4}}\right)\frac{1}{q_k^2}.\end{equation}
Next, we define
\[t(n)=  \frac{\varepsilon + \frac{1}{2\sqrt{b^2+4}}}{n}, \quad t'(n) = \frac{\varepsilon+\frac{q_k\delta_k}{2}}{n}.\]
We argue similarly to \eqref{step_3(1)} to show that
\begin{equation}\label{t_final}
    t'(n) \gs t(n) - \frac{1}{(b^2+4)^{3/2}}\frac{1}{nq_k^2}.
\end{equation}
Finally, combining \eqref{s_final} and \eqref{t_final} and assuming that $n$ is large enough such that
\mbox{$s(n)^2 - t(n)^2 \gs 1/2$}, we can deduce that

\begin{equation}
    \label{final_s_t}
    s'(n)^2 - t'(n)^2
    \gs \Big({s(n)^2-t(n)^2}\Big)\Bigg(1 - 4\Bigg(\Big(\frac{2}{(b^2+4)^{3/2}}+ \frac{1}{\sqrt{b^2+4}}\Big)\frac{1}{q_k^2}\Bigg)\Bigg).
\end{equation}
For $k$ sufficiently large, the error term in \eqref{final_s_t} is bounded from below by $1/2$, 
so we can apply Proposition \ref{exp_log_lemma} (iii) together with \eqref{middle_sized_2} to obtain

\begin{align}
    \label{final_step_3}
    \prod\limits_{n=N_0}^{\psi(q_k)} p_n(k,\varepsilon)  &\gs
    \Bigg(\prod\limits_{n=N_0}^{\psi(q_k)}g_n(\varepsilon)\Bigg) \Bigg(1 - \Big(a\frac{(\psi(q_k)+1)^2}{q_k^2} + 2b\frac{\psi(q_k)}{3q_k^2}\Big)
\Bigg)\times
\\&
\Bigg(1 - 8\Bigg(\left(\frac{2}{(b^2+4)^{3/2}}+ \frac{1}{\sqrt{b^2+4}}\right)\frac{\psi(q_k)}{q_k^2}\Bigg)\Bigg)
    = \Bigg(\prod\limits_{n=N_0}^{\psi(q_k)}g_n(\varepsilon)\Bigg)\Big(1 - \mathcal{O}\Big(\frac{\psi(q_k)^2}{q_k^2}\Big)\Big).
    \nonumber
\end{align}\vspace{3mm}\par

\subsection*{Proof of Step 5}
We will imitate a part of the argument from Step 4. However, here we cannot guarantee that the assumptions of Lemma 
\ref{better_sin} are fulfilled. So we use the weaker estimate 
\begin{equation}\label{sin_estim_easy}t^2(1 - 2\lvert t\rvert) \ls \sin^2(t) \ls t^2, \quad \lvert t \rvert < 1\end{equation}
on $x,y,z$ to obtain

\begin{align}\label{small_n_1}\nonumber&\prod_{n =1}^{N_0-1}\frac{f^2\Bigg(\frac{n}{q_k} - \Big(\Big\{\frac{nq_{k-1}}{q_k}\Big\}-\frac{1}{2}\Big)\delta_k\Bigg)}{f^2\big(\frac{n}{q_k}\big)}
 - \frac{f^2\left( \frac{2\varepsilon+q_k\delta_k}{2q_k}\right)}{f^2\big(\frac{n}{q_k}\big)}
\gs \\&\prod_{n =1}^{N_0-1}\frac{(x+y)^2 - z^2}{x^2}\prod_{n =1}^{N_0-1} \Bigg(1 + \mathcal{O}\Big(\frac{n^2}{q_k^2}\Big)\Bigg) =
\left(\prod_{n =1}^{N_0-1}\frac{(x+y)^2 - z^2}{x^2}\right)\left(1 + \mathcal{O}\big({q_k^{-2}}\big)\right)
\end{align}
where the implied constant depends on $N_0$.
By the same arguments as in Step 4, we can prove the inequalities \eqref{s_final} and \eqref{t_final}.
Since we cannot assume that $s(n)^2 - t(n)^2 \gs 1/2$, we now get from

\begin{align}
    \prod_{n=1}^{N_0-1} s'(n)^2 - t'(n)^2 \gs \prod_{n=1}^{N_0-1} \left(s(n)^2 - t(n)^2 + \mathcal{O}\big({q_k^{-2}}\big)\right)
    = \left(\prod_{n=1}^{N_0-1} \big(s(n)^2 - t(n)^2\big)\right) + \mathcal{O}\big({q_k^{-2}}\big)
\end{align}
an \textit{additive} error of order $\mathcal{O}\big({q_k^{-2}}\big)$.

\subsection*{Proof of Step 6}
We know that there exists a constant $c_4 > 0$ such that
\[c_4 = \max_{k \gs 1} q_k\delta_k + \max_{\varepsilon \in I}\varepsilon.\]
 Using \eqref{sin_estim_easy},
 we get
 
 \begin{align}\label{step_5_final}f\Big(\delta_k + \frac{\varepsilon}{q_k}\Big)q_k 
 &\gs
 2\pi \Big\lvert q_k\delta_k + \varepsilon\Big\rvert \left(1 - \pi^2(q_k\delta_k + \varepsilon)^2q_k^{-2}\right)\\&\gs 
 2\pi \Big\lvert \frac{1}{\sqrt{b^2+4}} + \varepsilon\Big\rvert\left(1 - \pi^2c_4^2q_k^{-2}\right)
 - 2\pi \Big\lvert q_k\delta_k - \frac{1}{\sqrt{b^2+4}}\Big\rvert
 \left(1 - \pi^2c_4^2q_k^{-2}\right).
 \end{align}
  Arguing as in \eqref{step_3(1)}, we have
 \[\left\lvert q_k\delta_k - \frac{1}{\sqrt{b^2+4}}\right\rvert \ls \frac{2}{b^2 +4}q_k^{-2} = \mathcal{O}(q_k^{-2}),\]
 hence
 
 \[f\Big(\delta_k + \frac{\varepsilon}{q_k}\Big)q_k \gs 
  2\pi \Big\lvert \frac{1}{\sqrt{b^2+4}} + \varepsilon\Big\rvert\left(1 + \mathcal{O}(q_k^{-2})\right)
 + \mathcal{O}(q_k^{-2})
 \]
 follows.
 This finishes the proof of Lemma \ref{quantitative_convergence}.

\subsection*{Proof of Lemma \ref{truncation_lem}}
The proof of Lemma \ref{truncation_lem} works in a similar fashion to Step 3 of Lemma \ref{quantitative_convergence}, although this time we need a \textit{lower} bound on $\prod\limits_{n = T+1}^{\infty}g_n(\beta,\varepsilon)$.
If $T$ is sufficiently large, we can remove the absolute values in the definition of $g_n(\varepsilon)$ and obtain
\begin{align}
    g_n(\varepsilon) 
& = 1 - 2\ell(n)
+ \left(\frac{1}{b^2 + 4}\big(\{n\beta\} - \tfrac{1}{2}\big)^2 - \Big(\varepsilon + \frac{1}{2\sqrt{b^2+4}}\Big)^2\right)\frac{1}{n^2}
\\& \gs 1 - \Big(2\ell(n) + \frac{M}{n^2}\Big)
\end{align}
where $\ell(n)$ is defined as in \eqref{ell_n}
and 
$M = \big(\frac{1}{4(b^2+4)} + \max\limits_{\varepsilon \in I}
\varepsilon + \frac{1}{2\sqrt{b^2+4}}\big)^2$.
Using \eqref{sum_k_bound}, we get for $b >1$

\begin{equation}\label{lem5_bound}\Big\lvert\sum_{n = T+1}^{\infty} 2\ell(n) + \frac{M}{n^2}\Big\rvert 
\ls \frac{2b}{\sqrt{b^2 +4}\log(b)}\frac{\log(T +1)}{T} + \frac{2 + M}{T}.
\end{equation}
We see that for $T$ sufficiently large, 
\eqref{lem5_bound} can be bounded by $1/2$, hence applying \mbox{$\text{Proposition \ref{exp_log_lemma} (ii)}$} leads to

\begin{align}\label{final_truncation}
    \prod\limits_{n = T+1}^{\infty}g_n(\beta,\varepsilon)
    \gs 1 - \Bigg(\frac{3b}{\sqrt{b^2 +4}\log(b)}\frac{\log(T +1)}{T} + \frac{6 + 3M}{2T}\Bigg).
\end{align}
Again, for $b =1$, we can replace $\log(b)$ in the estimate by $\log(3/2)$.

\subsection*{ Proof of Corollary \ref{quantitative_convergence_cor}}

We can explicitly compute the constants $c_1,\ldots, c_5$ used in the proof of Lemma \ref{quantitative_convergence} for the Golden Ratio respectively the quadratic irrational $\beta(5)$, and specify all the errors from Steps 1\,--\,6. It turns out that $T = 100.000$ and $k=24$ for the Golden Ratio respectively $k=10$ for $\beta(5)$ are large enough to fulfill
all assumptions on $k,T$ to be ``sufficiently large'' that are made in the proofs of Lemmas \ref{quantitative_convergence} and \ref{truncation_lem}. 
We see that the errors from Steps $1,5,6$ are in $\mathcal{O}(q_k^{-2})$,
and since all constants are small and since we can choose in both cases $N_0=2$, 
we can bound both the additive and the multiplicative error by $0.0001$ in the Golden Ratio case and
by $0.00001$ in the case $b = 5$.
To not repeat all arguments from above again, we 
will only address Steps $2,3,4$ and consider the main error terms there.
In the case of the Golden Ratio, we can bound $c_1 \ls 0.8, c_2 \ls 0.5, c_3 \ls 4$, and thus we get for Step 2 an error bounded by $\frac{\frac{3}{2}(k+4) + 17}{\psi(q_k) -1} \ls 0.02.$
Step 3 leads to an error bounded by
$3.5\frac{\log(\psi(q_k) +1)}{\psi(q_k)} + \frac{3}{\psi(q_k)} \ls 0.01$, Step 4 to
\mbox{$\frac{12(\psi(q_k)+1)^2}{q_k^2} + \frac{16\psi(q_k)}{3q_k^2} \ls 0.01$}.
Adding these errors, we obtain the result for the Golden Ratio. Similarly, we can treat the case $b = 5$:
Here we can bound  $c_1 \ls 1.03, c_2 \ls 0.2, c_3 \ls 5$ which leads for Step 2 to an error bounded by
$\frac{14(k+4)}{5\psi(q_k)} + \frac{4/25 + 20}{\psi(q_k) -1} \ls 0.00045.$
For Step 3, we can bound the error by $\frac{15}{\sqrt{29}\log(5)}\frac{\log(\psi(q_k) +1)}{\psi(q_k)} + \frac{3}{\psi(q_k)}\ls 0.0004$, for Step 4 by
$\frac{9(\psi(q_k)+1)^2}{q_k^2} + \frac{18\psi(q_k)}{q_k^2} \ls 0.0011$.
Again, the result follows immediately from adding these errors.

\section{Acknowledgements}
The author would like to thank Christoph Aistleitner for helpful discussions and for suggesting this direction of research.

\end{document}